\documentclass[reqno,12pt]{amsart}
\usepackage[T1]{fontenc}
\usepackage{lmodern}
\usepackage{amsmath,amssymb,amsthm,mathrsfs}
\usepackage[margin=1in]{geometry}
\usepackage[expansion=false]{microtype}
\usepackage{graphicx,booktabs}
\usepackage{placeins}
\usepackage{hyperref}
\hypersetup{hidelinks}
\newtheorem{theorem}{Theorem}[section]
\newtheorem{lemma}[theorem]{Lemma}
\newtheorem{proposition}[theorem]{Proposition}
\newtheorem{corollary}[theorem]{Corollary}
\theoremstyle{definition}

\theoremstyle{remark}
\newtheorem{remark}[theorem]{Remark}
\numberwithin{equation}{section}
\newcommand{\R}{\mathbb R}
\newcommand{\C}{\mathbb C}
\newcommand{\Sph}{\mathbb S^2}
\newcommand{\supp}{\operatorname{supp}}
\newcommand{\eps}{\varepsilon}
\newcommand{\cR}{\mathcal R}

\newcommand{\norm}[1]{\left\lVert#1\right\rVert}

\newcommand{\ip}[2]{\left\langle#1,#2\right\rangle}
\newcommand{\dd}{\,\mathrm d}
\allowdisplaybreaks[4]
\makeatletter
\renewcommand{\@captionfont}{\normalfont\small}
\makeatother
\title[Inverse potential scattering]
{Increasing stability for the inverse potential problem of the fractional Helmholtz equation}
\author[D. Qiu]{Dong Qiu}
\address{School of Mathematical Sciences, Zhejiang University, Hangzhou 310058, China}
\email{qiudong@zju.edu.cn}
\author[X. Xu]{Xiang Xu}
\address{School of Mathematical Sciences, and Center for Interdisciplinary Applied Mathematics, Zhejiang University, Hangzhou 310058, China}
\email{xxu@zju.edu.cn}
\author[Y. Zhao]{Yue Zhao}
\address{School of Mathematics and Statistics, and Key Lab NAA-MOE, Central China Normal University, Wuhan 430079, China}
\email{zhaoy@ccnu.edu.cn}
\date{}

\begin{document}
\begin{abstract}
This paper is concerned with the inverse potential scattering problem for the three-dimensional fractional Helmholtz equation.
For every $0<\alpha<1$, we establish an increasing stability estimate for determining the potential function from multifrequency exterior measurements of the scattered field.
The stability estimate consists of a Lipschitz data discrepancy and a logarithmic remainder that decreases when the upper measured wavenumber increases.
The key ingredient in the analysis is to construct appropriate geometric optics solutions, which relate the X-ray transform of the potential to the measured fields.
As a consequence, the stability is derived by applying analytic continuation.
Two-dimensional numerical examples with noisy near-field data illustrate the benefit of enlarging the measured frequency band for three representative fractional orders.
\end{abstract}
\maketitle

\section{Introduction}\label{sec:intro}

As the basis of fractional quantum mechanics, the fractional Schr\"odinger equation was originally established by Laskin via the Feynman path-integral formalism for L$\acute{\rm e}$vy-flight particles \cite{Laskin1, Laskin2}.
Concurrently, its associated fractional Laplacian has emerged as a powerful effective model capable of capturing intrinsic non-locality, power-law dispersion, and anomalous attenuation that arise in a variety of physical systems, including quantum optics \cite{S1, S2}, geophysical electromagnetics \cite{sisc}, and electrodynamics \cite{TT}.
Driven by the important applications, the fractional Schr\"odinger equation has been extensively investigated in both mathematical and physical regimes, see e.g., \cite{LL, Liu, naber, prl} and the references cited therein.

We consider the equation
\begin{equation}\label{eq:model}
(-\Delta)^\alpha u-k^{2\alpha}u+Vu=0 \quad\text{in }\R^3,\qquad 0<\alpha<1,
\end{equation}
where $k>0$ is the wavenumber and $V$ is a real-valued, compactly supported potential independent of $k$.
We use the Fourier transform convention for the fractional Laplacian operator
\[
\widehat f(\xi)=\int_{\R^3}e^{-ix\cdot\xi}f(x)\dd x, \qquad \widehat{(-\Delta)^\alpha f}(\xi) =|\xi|^{2\alpha}\widehat f(\xi).
\]
For the incident plane wave $u^{\rm inc}=e^{ikx\cdot d}$, $d\in\Sph$, write the total field as $u=u^{\rm inc}+u^s$.
The scattered field is outgoing and satisfies
\begin{equation}\label{eq:radiation}
u^s(\rho\theta)=O(\rho^{-1}),\qquad \rho(\partial_\rho-ik)u^s(\rho\theta)\longrightarrow0 \quad\text{as }\rho\to\infty,
\end{equation}
uniformly in $\theta\in\Sph$, where $\rho = |x|$ and $\theta\in\mathbb S^2$.
The direct problem is to determine the scattered field $u^s$ for a given potential and incident wave, while the inverse problem is to reconstruct $V$ from exterior measurements of $u^s$.

For the direct scattering problems, a great deal of mathematical and numerical results are available involving the fractional Laplacian, in particular the fractional Schr\"odinger equation.
In quantum optics, kinetic equations for two-photon states in random media were derived in \cite{S2}.
Collective spontaneous emission and kinetic equations for one-photon states were studied in \cite{S1}.
Complementary to the results in mathematical theory of the scattering resonances for the classical Schr\"odinger operator, meromorphic continuation for the resolvent of the fractional Schr\"odinger operators was investigated in \cite{LYZ}, which yields a logarithmic analytic region and resolvent estimates.
Numerical methods for fractional Helmholtz equations arising in electromagnetic modeling were developed in \cite{bl, sisc}.

Inverse problems for the fractional Schr\"odinger equation have received considerable attention in recent years.
Much of the existing work concerns the fractional Calder\'on problem, which seeks to recover a potential from a nonlocal exterior Dirichlet-to-Neumann map.
Ghosh, Salo, and Uhlmann proved global uniqueness for the fractional Calder\'on problem using unique continuation and Runge approximation \cite{GSU}.
Further uniqueness and stability results were obtained in \cite{grsu, ruland}.
Covi, de Hoop, and Salo constructed geometrical optics solutions for the fractional Helmholtz equation and applied them to inverse potential problems \cite{CDS}.
Compared with the well-studied fractional Calder${\rm \acute o}$n problems, little is known for the uniqueness and stability for the inverse scattering problems modeled by the fractional Schr\"odinger equation. 
Uniqueness for inverse potential scattering from multifrequency far-field data was proved in \cite{DGM}.
The principal symbols of random-source covariance and relation operators were recovered almost surely from multifrequency far-field data generated by a single realization in \cite{LL}.
Increasing stability for the inverse source problem of the stochastic fractional Helmholtz equation was established in \cite{QXZ}.
Bayesian approaches to inverse scattering for fractional Helmholtz models were studied in~\cite{Jia}, and the factorization method for support reconstruction was developed in~\cite{ZilberbergSupport}. For related reconstruction methods and Bayesian analysis for the fractional Calder\'on problem, see also~\cite{DRR, KNR, LiSingle, LWF}.
The present work concerns stability for recovering a deterministic potential from exterior scattered fields.

We formulate the main result of this paper. Fix $0<r<R$ and $F>0$, and introduce the a priori class
\begin{equation}\label{eq:class}
\mathcal V_{r,F}^4 =\{V\in C_c^4(B_r;\R): \norm V_{C^4(\R^3)}\le F\}.
\end{equation}
Here $B_r=\{x:|x|<r\}$ and $\norm V_{C^4}=\max_{|\nu|\le4}\norm{\partial^\nu V}_\infty$.
The observations are made either on an enclosing sphere or in an exterior annulus.
Set
\begin{equation}\label{eq:geometry}
\Gamma=\partial B_R,\qquad A=B_R\setminus\overline{B_r},\qquad \mathcal O\in\{\Gamma,A\},\qquad Y=L^2(\mathcal O\times\Sph).
\end{equation}
The observation measure is surface measure on $\Gamma$ or volume measure on $A$, together with surface measure in $d$.
For two potentials let $u_j^s=u_{V_j}^s$, $j=1,2$, and put $\tau=(2\alpha-1)_+=\max(2\alpha-1,0)$.

The following increasing stability for the inverse potential scattering problem is the main result of this paper. It shows that the stability 
will improve as the upper bound of the bandwidth of the utilized frequencies increases.
\begin{theorem}\label{thm:main-stability}
Let $0<\alpha<1$ and $V_1,V_2\in\mathcal V_{r,F}^4$.
There are constants $K_0,C>0$ such that the following holds for $K\ge K_0+1$.
Define
\begin{equation}\label{eq:data-discrepancy}
\eps_w=\sup_{k\in[K_0,K]} k^\tau \norm{u_1^s(k)-u_2^s(k)}_Y.
\end{equation}
If $0<\eps_w<e^{-1}$, set $E=|\log\eps_w|$ and $\rho=\min(8\alpha,4)$.
Then the following stability estimate holds
\begin{equation}\label{eq:main-stability}
\norm{V_1-V_2}_{H^{-1/2}(\R^3)}^2 \le C\left(\eps_w^2+ \left(K+\sqrt{\frac{E}{\log(2+E)}}\right)^{-\rho}\right).
\end{equation}
The constants are independent of $K$ and $\eps_w$. As a consequence, an increasing stability estimate also holds for the $L^2$ norm of the potential
by an argument of interpolation.  
\end{theorem}

In this work, we establish high-frequency well-posedness of the direct scattering problem for every \(0<\alpha<1\), together with uniform resolvent estimates and holomorphic continuation in a fixed complex strip.
For the inverse problem, our approach employs the construction of the geometric optics solutions, which relate the X-ray transform of the potential
to the exterior measurements. To do this, we construct an exact real phase and a positive leading amplitude for the geometric optics solutions.
The fractional order determines the size of the leading phase correction, which is proportional to $k^{1-2\alpha}$ times a half-ray integral of the potential.
This correction tends to zero for $\alpha>1/2$, is of order one for $\alpha=1/2$, and may grow for $\alpha<1/2$.
In the last case, its exponential may wind around the unit circle.
To deal with this issue, 
we use the integer combinations of phases to cancel the second phase coefficient, and the frequency averaging method to control the phase ambiguity 
for $\alpha<1/2$. As a consequence, the Fourier transform of the potential can be controlled by the averaged frequency data with a high-frequency remainder term.
Then the stability follows from a quantitative analytic continuation argument as in Bao, Zhai, and Zhao~\cite{BZZ}.
Finally, we implement a regularized multifrequency reconstruction method for the two-dimensional model.
Numerical examples for different fractional orders illustrate the bandwidth dependence of the reconstruction accuracy.

This paper is organized as follows. Section~\ref{sec:direct} establishes direct scattering and the resolvent estimates.
Section~\ref{sec:inverse} constructs the geometric optics approximation and proves its error estimates in the volume and on the receiver sphere.
Section~\ref{sec:stability} relates the phases to ray integrals and proves Theorem~\ref{thm:main-stability} by frequency combinations and analytic continuation.
Section~\ref{sec:dense-numerics} presents the reconstruction method and numerical examples.
Section~\ref{sec:conclusion} concludes the paper, and the appendices provide some technical proofs.

Throughout the paper, $C>0$ denotes a generic constant that may change from line to line.
Unless otherwise stated, it depends only on the fixed fractional order, geometry, and a priori bound.
We also denote $\langle x\rangle=(1+|x|^2)^{1/2}$, $D=-i\nabla$, $t_-=\max(-t,0)$, and $\ip f g=\int f\overline g$.

\section{Direct scattering problem}\label{sec:direct}

We establish high-frequency well-posedness and the resolvent bounds needed for the field approximation and analytic continuation.
As we are focused on the inverse problem, proofs for some results in this section are removed to the appendix, which are related closely to the scattering theory.
We refer to the monograph \cite{DZ} for a comprehensive introduction to the scattering theory.

\subsection{Green function and well-posedness}

Denote the outgoing free resolvents by
\[
\cR_{\alpha,0}(k)=\bigl((-\Delta)^\alpha-k^{2\alpha}-i0\bigr)^{-1}, \qquad R_0(k)=(-\Delta-k^2-i0)^{-1}.
\]
In the right half-plane, complex powers use $-\pi/2<\arg z<\pi/2$.
For the free equation, we work with tempered distributions whose Fourier transforms are locally integrable near the origin.
The fractional multiplier is defined by multiplication there and by its smooth action on distributions away from the origin.
Outgoing source-generated fields additionally satisfy \eqref{eq:radiation}.
The free resolvent estimates used below were established in \cite[Proposition 2.7]{QXZ}.
We recall them in the form needed here.
The following decomposition and estimates are proved in Appendix~\ref{app:free}.

\begin{lemma}\label{lem:free}
Let $0<\alpha<1$.
The outgoing free Green function in three dimensions is
\begin{equation}\label{eq:green}
G_{\alpha,k}(x)=\frac{k^{2-2\alpha}}{\alpha} \frac{e^{ik|x|}}{4\pi|x|}+J_{\alpha,k}(|x|),
\end{equation}
where
\begin{equation}\label{eq:J-kernel}
J_{\alpha,k}(a)=\frac{\sin(\pi\alpha)}{4\pi^2a} \int_0^\infty \frac{\gamma^\alpha e^{-a\sqrt\gamma}} {\gamma^{2\alpha}-2\gamma^\alpha k^{2\alpha}\cos(\pi\alpha)+k^{4\alpha}} \dd\gamma.
\end{equation}
For smooth compactly supported cutoffs $\chi_1,\chi_2$, $\chi_1\cR_{\alpha,0}(z)\chi_2$ is holomorphic from $L^2$ to $H^{2\alpha}$ on $\Re z>0$.
For $|\arg z|\le\theta_0<\pi/2$, $|z|\ge1$, and $0\le l\le2\alpha$,
\begin{equation}\label{eq:free-bound}
\norm{\chi_1\cR_{\alpha,0}(z)\chi_2}_{L^2\to H^l} \le C|z|^{1+l-2\alpha}e^{T(\Im z)_-},
\end{equation}
where $T>0$ depends on the cutoffs.
If $J_\alpha(z)$ denotes the convolution operator associated with the second term, then
\begin{equation}\label{eq:J-bound}
\norm{J_\alpha(z)}_{L^2\to H^l}\le C|z|^{l-2\alpha}.
\end{equation}
For each fixed $k>0$, $G_{\alpha,k}*f$ is the unique outgoing solution of the free equation with a compactly supported $L^2$ source $f$, and belongs to $H^{2\alpha}_{\mathrm{loc}}(\R^3)$.
\end{lemma}

The Lippmann--Schwinger equation associated with \eqref{eq:model} is
\begin{equation}\label{eq:LS}
u=e^{ikx\cdot d}-\cR_{\alpha,0}(k)(Vu),\qquad u^s=-\cR_{\alpha,0}(k)(Vu).
\end{equation}
By Lemma~\ref{lem:free}, this integral equation is equivalent to \eqref{eq:model} with the radiation condition \eqref{eq:radiation} in the distributional class specified above.

We denote the outgoing perturbed resolvent by
\[
\cR_{\alpha,V}(z)=\bigl((-\Delta)^\alpha+V-z^{2\alpha}\bigr)^{-1}_{\rm out}.
\]
For $z$ in the first quadrant, it agrees with the spectral resolvent.
Elsewhere, the localized resolvent is defined by outgoing meromorphic continuation.

\begin{theorem}\label{thm:direct}
Let $0<\alpha<1$ and $V\in L^\infty_c(\R^3;\R)$.
There is a discrete set $\Sigma_V\subset(0,\infty)$ such that, for $k\notin\Sigma_V$, \eqref{eq:LS} has a unique solution with $u^s\in H^{2\alpha}_{\mathrm{loc}}(\R^3)$.
The scattered field depends continuously in $H^{2\alpha}_{\mathrm{loc}}$ on the $L^2$ restriction of the incident field to a fixed neighborhood of $\supp V$.
The localized outgoing resolvent is meromorphic on $\Re z>0$.

If in addition $V\in C_c^1(\R^3;\R)$, then
\begin{equation}\label{eq:threshold}
\Sigma_V\subset\{k>0:k^{2\alpha}\le E_*(V)\},
\end{equation}
where
\begin{equation}\label{eq:intro-threshold}
E_*(V)=\max\left(0,\sup_{x\in\R^3} \left(V(x)+\frac{x\cdot\nabla V(x)}{2\alpha}\right)\right).
\end{equation}
In particular, the direct problem is well posed at every sufficiently large frequency for every $0<\alpha<1$.
\end{theorem}

The proof for this theorem provided in Appendix~\ref{app:direct} uses analytic Fredholm theory to construct the outgoing solution.
At an exceptional real frequency, the corresponding homogeneous solution is an eigenfunction; a scaling identity gives \eqref{eq:threshold}.

\begin{remark}
For potentials supported in $B_r$ with a common $C^1$ bound, the threshold is uniform, since
\[
E_*(V)\le\norm V_\infty+\frac r{2\alpha}\norm{\nabla V}_\infty.
\]

\end{remark}

\subsection{High-frequency resolvent estimates}\label{sec:resolvent}

For the remainder of the paper, the potentials belong to $\mathcal V_{r,F}^4$.
Put $\mu=2\alpha$.
The following weighted estimate is proved by a positive commutator argument in Appendix~\ref{app:real}.

\begin{lemma}\label{lem:real-resolvent}
Let $0<\mu<2$ and $V\in\mathcal V_{r,F}^4$.
There are $k_0,C>0$ such that the outgoing resolvent satisfies
\begin{equation}\label{eq:weighted-real}
\norm{\langle x\rangle^{-1}\cR_{\alpha,V}(k) \langle x\rangle^{-1}}_{L^2\to L^2} \le Ck^{1-\mu},\qquad k\ge k_0.
\end{equation}
The same estimate holds with two fixed compactly supported cutoffs in place of the spatial weights.
\end{lemma}

The free-resolvent decomposition gives a fixed strip of holomorphy for the localized perturbed resolvent; see Appendix~\ref{app:strip}.

\begin{proposition}\label{prop:complex-strip}
Let $0<\mu<2$ and fix compactly supported smooth cutoffs $\eta_1,\eta_2$.
There exist $k_0,b,C>0$ such that, uniformly for $V\in\mathcal V_{r,F}^4$, the localized outgoing resolvent $\eta_1\cR_{\alpha,V}(z)\eta_2$ is holomorphic in
\begin{equation}\label{eq:strip}
\mathcal S_{k_0,b}=\{z\in\C:\Re z>k_0,\ |\Im z|<b\}.
\end{equation}
Moreover,
\begin{equation}\label{eq:complex-bound}
\norm{\eta_1\cR_{\alpha,V}(z)\eta_2}_{L^2\to L^2} \le C(\Re z)^{1-\mu},\qquad z\in\mathcal S_{k_0,b}.
\end{equation}
\end{proposition}

\begin{corollary}\label{cor:data-analyticity}
Let $Y$ be defined by \eqref{eq:geometry}.
Uniformly for $V\in\mathcal V_{r,F}^4$, the map $z\mapsto u_V^s(\cdot,\cdot,z)$ is $Y$-valued holomorphic in a strip of the form \eqref{eq:strip}, and
\begin{equation}\label{eq:polynomial-data}
\norm{u_V^s(\cdot,\cdot,z)}_Y\le C(1+|z|)^3.
\end{equation}
The total field and its modulation by $e^{-izx\cdot d}$ have the same properties.
\end{corollary}

The volume estimate follows from Proposition~\ref{prop:complex-strip}.
For spherical observations, the separation of the receivers from the potential support gives a holomorphic integral kernel; see Appendix~\ref{app:data-analyticity}.

\section{Geometric optics approximation}\label{sec:inverse}

In this section, we investigate the construction of geometric optics solutions to the fractional Helmholtz equation, which is useful in the subsequent analysis of the inverse problem.
To do this,
we construct a real phase and a positive leading amplitude, then estimate the difference between the resulting approximation and the outgoing field.
For background on semiclassical analysis and geometric optics, we refer to~\cite{Zworski}.
Throughout this section, $\mu=2\alpha$ and $Y$ is defined by \eqref{eq:geometry}.

For a compactly supported function $f$, define its half-ray and X-ray transforms by
\begin{equation}\label{eq:beam-xray}
Bf(x,d)=\int_{-\infty}^0f(x+sd)\dd s, \qquad Xf(x,d)=\int_\R f(x+sd)\dd s.
\end{equation}
They satisfy $d\cdot\nabla Bf=f$ and $Xf(x,d)=Bf(x,d)+Bf(x,-d)$.
Put $h=k^{-1}$, $\delta=h^\mu$, and
\[
P_{h,V}=|hD|^\mu-1+\delta V, \qquad R_{h,V}=(P_{h,V}-i0)^{-1}.
\]
Rescaling Lemma~\ref{lem:real-resolvent} gives
\begin{equation}\label{eq:real-input}
\norm{\langle x\rangle^{-1}R_{h,V} \langle x\rangle^{-1}}_{L^2\to L^2}\le Ch^{-1}.
\end{equation}

\subsection{The phase and leading amplitude}\label{sec:phase-construction}

The incoming Hamilton flow defines the phase through its action and a positive leading amplitude through the transport equation.

\begin{lemma}\label{lem:phase-construction}
Let $0<\mu<2$ and $V\in\mathcal V_{r,F}^4$.
There are $\delta_0,C>0$, depending only on $\mu,r,R,F$, such that for $|\delta|\le\delta_0$ there are an incoming phase $S_{\delta,V,d}$ and a leading amplitude $a_{0,\delta,V,d}$ with the following properties.
The phase is real, and both $S_{\delta,V,d}$ and its spatial gradient have uniformly bounded derivatives through order three jointly in $(x,d,\delta)$ on the fixed construction sets.
The leading amplitude is real and positive.
On a fixed neighborhood of $\overline{B_R}$,
\begin{equation}\label{eq:eikonal}
|\nabla S_{\delta,V,d}|^\mu=1-\delta V,\qquad \nabla S_{\delta,V,d}=d+O(\delta),\qquad a_{0,\delta,V,d}=1+O(\delta)\ge\tfrac34.
\end{equation}
Moreover,
\begin{equation}\label{eq:phase-expansion}
\begin{split}
S_{\delta,V,d}-x\cdot d&=-\frac\delta\mu BV(x,d)+O(\delta^2),\\
\partial_\delta S_{\delta,V,d}&=-\frac1\mu BV(x,d)+O(\delta).
\end{split}
\end{equation}
These remainders and the amplitude remainder are uniform in the $C^1$ norm in the spatial and incident-direction variables.
The phase and leading amplitude equal $x\cdot d$ and one in the incoming region $x\cdot d<-r$.
\end{lemma}

\begin{proof}
Treat $\delta$ as an independent parameter and prescribe $X(0)=y-L_0d$, $\Xi(0)=d$ on an incoming section, with $y\in d^\perp$ and $L_0>r$ fixed.
The Hamilton equations are
\[
\dot X=\mu|\Xi|^{\mu-2}\Xi,\qquad \dot\Xi=-\delta\nabla V(X).
\]
The Hamiltonian $H_\delta=|\Xi|^\mu+\delta V(X)$ equals one along every trajectory.
For sufficiently small $|\delta|$, the momentum remains in $1/2<|\Xi|<2$ on the fixed time interval.
The vector field and its flow have uniform joint $C^3$ bounds.
The projection $(y,t)\mapsto X(t;y,d,\delta)$ is a small $C^1$ perturbation of the affine map $y-L_0d+\mu td$.
It has an inverse on the fixed construction region, after enlarging the parameter sets by a fixed margin.
The inverse time and incoming label $y=y(x,d,\delta)$ are jointly $C^3$.
Finitely many transverse frames cover the incident sphere.

Define the reduced action by
\[
S_{\delta,d}(x)=-L_0+\int_0^{t(x,d,\delta)}\Xi(s)\cdot\dot X(s)\dd s.
\]
Differentiating along a family of trajectories at fixed $\delta$ and integrating by parts give
\[
\dd S=\Xi(t)\cdot\dd x-d\cdot\dd X(0)+\int_0^t\dd H_\delta(X(s),\Xi(s))\dd s.
\]
The terminal-time variation is included in $\dd x$; the last integral vanishes because every trajectory has energy one.
For fixed $d$, the initial variation satisfies $d\cdot\dd y=0$.
For fixed $x$ and a tangential variation of $d$, use $\dd X(0)=\dd y-L_0\dd d$ and $d\cdot\dd y=-y\cdot\dd d$.
Thus
\begin{equation}\label{eq:momentum-regularity}
\nabla_xS=\Xi,\qquad \nabla_{\Sph,d}S=y,\qquad S,\Xi,y\in C^3_{x,d,\delta}.
\end{equation}
The intrinsic spherical gradient $y$ lies in $d^\perp$.
The identities
\[
\Xi(t)-d=-\delta\int_0^t\nabla V(X(s))\dd s,\qquad
S-x\cdot d=\int_0^{t(x,d,\delta)}(\Xi(s)-d)\cdot\dot X(s)\dd s
\]
give $\Xi-d=O_{C^3_{x,d}}(\delta)$ and $S-x\cdot d=O_{C^3_{x,d}}(\delta)$ at fixed $\delta$.

We also need the mixed bound
\begin{equation}\label{eq:phase-third-mixed}
\sup_{|\delta|\le\delta_0}\norm{\partial_\delta^3S}_{C^1_{x,d}}\le C.
\end{equation}
The joint $C^3$ bound gives a continuous third parameter derivative.
For a fixed spatial path $\gamma$, integrate the first identity in \eqref{eq:momentum-regularity} along $\gamma$ and differentiate three times in $\delta$.
Since $\Xi$ is jointly $C^3$, this gives $\nabla_x\partial_\delta^3S=\partial_\delta^3\Xi$.
The same argument along paths in a spherical chart, using the incoming label, gives $\nabla_{\Sph,d}\partial_\delta^3S=\partial_\delta^3y$.
Both right-hand sides are continuous and uniformly bounded, proving \eqref{eq:phase-third-mixed}.
Differentiating the eikonal equation at zero yields $\partial_\delta S_0=-BV/\mu$, and Taylor's formula proves \eqref{eq:phase-expansion}.

Choose a smooth real radial symbol $p_{\rm sm}$ which is zero near the origin and equals $|\xi|^\mu$ for $|\xi|\ge1/4$.
Put $v=\nabla_\xi p_{\rm sm}(\Xi)$ and solve
\begin{equation}\label{eq:transport-recursion}
v\cdot\nabla a_0+\tfrac12(\operatorname{div}v)a_0=0,\qquad a_0=1\quad\text{upstream}.
\end{equation}
Here $v$ is jointly $C^3$ and $\operatorname{div}v$ is jointly $C^2$.
The characteristic formula
\[
a_0(x,d,\delta)=\exp\left(-\tfrac12\int_0^{t(x,d,\delta)}\operatorname{div}v(X(s),d,\delta)\dd s\right)
\]
gives a joint $C^2$ bound and, after reducing $\delta_0$, the following estimates at fixed $\delta$:
\begin{equation}\label{eq:derivative-budget}
\Xi-d=O_{C^3_{x,d}}(\delta),\qquad a_0-1=O_{C^2_{x,d}}(\delta),\qquad a_0\ge\tfrac34.\qedhere
\end{equation}
\end{proof}

\subsection{Volume approximation}\label{sec:finite-wkb}

We use the following free propagation estimate for sources localized near the incident direction.
Its proof is given in Appendix~\ref{app:directed-trace}.

\begin{lemma}\label{lem:directed-free-trace}
Let $0<\mu<2$, $0<a<R$, and $T>0$ be fixed, and let $d\in\Sph$.
Write $x=p+td$, where $p\in d^\perp$, and set
\[
K_d=\{p+td:|p|\le a,\ |t|\le T\},\qquad \Gamma=\partial B_R.
\]
There is a fixed $\varepsilon>0$ with the following property.
Let $q_d\in C_c^\infty$ be supported in $|\xi-d|<\varepsilon$, with uniformly bounded symbol seminorms under rotations taking $d$ to a fixed axis.
Then, for $0<h\le1$ and $\operatorname{supp}f\subset K_d$,
\begin{equation}\label{eq:microlocal-sphere}
\|R_{h,0}q_d(hD)f\|_{L^2(\Gamma)}\le Ch^{-1}\|f\|_2,
\qquad R_{h,0}=(|hD|^\mu-1-i0)^{-1}.
\end{equation}
The source may meet $\Gamma$.
Moreover, if $K_D$ and $K_O$ are fixed bounded source and observation sets with
\[
y\cdot d\ge T_+\quad(y\in K_D),\qquad
x\cdot d\le T_-\quad(x\in K_O),\qquad T_-<T_+,
\]
then, for every fixed integer $N$,
\begin{equation}\label{eq:backward-cone}
\|R_{h,0}q_d(hD)f\|_{L^\infty(K_O)}\le C_Nh^N\|f\|_2,
\qquad \operatorname{supp}f\subset K_D.
\end{equation}
All constants are uniform in $d$.
\end{lemma}

The total field $u_{V,h,d}=u_V(\cdot,d,h^{-1})$ admits the following approximation.

\begin{lemma}\label{lem:wkb-l2}
Let $0<\mu<2$ and $V\in\mathcal V_{r,F}^4$.
For the phase and amplitude in Lemma~\ref{lem:phase-construction}, there are $h_0,C>0$ such that the total field satisfies
\begin{equation}\label{eq:wkb-l2}
u_{V,h,d}=e^{iS_{\delta,V,d}/h} (a_{0,\delta,V,d}+\rho_{V,h,d}),\qquad \sup_{d\in\Sph}\norm{\rho_{V,h,d}}_{L^2(B_R)} \le Ch^{1+\mu},\quad 0<h\le h_0.
\end{equation}
\end{lemma}

\begin{proof}
Write $u=u_{V,h,d}$ and abbreviate the phase and amplitude of Lemma~\ref{lem:phase-construction} by $S$ and $a_0$.
Extend $S-x\cdot d$ and $a_0-1$ by fixed cutoffs outside a larger construction slab, preserving them near the observation and cutoff sets.
For $S^{\rm ext}=x\cdot d+\chi(S-x\cdot d)$, the formula
\[
\nabla S^{\rm ext}=d+\chi(\Xi-d)+(\nabla\chi)(S-x\cdot d)
\]
retains the fixed-parameter $C^3$ momentum bound.
The corresponding extension $1+\chi(a_0-1)$ retains the $C^2$ amplitude bound and positivity.
We continue to write $S,a_0$ for these extensions.

\medskip
\noindent\textit{The directional decomposition.}
Choose fixed numbers
\[
r<r_0<r_1<r_2<R.
\]
With a sufficiently small fixed upper bound on $\delta$, trajectories at transverse distance at least $r_0$ from the incident axis do not reach $\supp V$ in the construction slab, since their displacement from free trajectories is uniformly $O(\delta)$ on the fixed time interval.
The extended function
\[
W_h=e^{iS/h}a_0-e^{ix\cdot d/h}
\]
therefore has support in a fixed bounded tube with $|x-(x\cdot d)d|\le r_0$.
The extension cutoffs preserve this support property.

Let $q_d\in C_c^\infty(\mathbb R^3)$ equal one in a neighborhood of $d$ and have support in a sufficiently small fixed ball about $d$.
Its support is disjoint from $|\xi|\le1/4$.
Decrease the admissible upper bound on $\delta$ so that $q_d=1$ on a fixed neighborhood of every extended phase gradient.
Put
\[
Q_h=q_d(hD),\qquad Q_h^\perp=I-Q_h.
\]
Choose this radius small relative to $r_2-r_1$ and the longitudinal size of the construction region.
All these choices are independent of $h$ and of the particular potential.

Choose $\chi_*\in C_c^\infty$ equal to one near $\overline{B_R}$, with downstream derivatives in a compact subslab of $2R<x\cdot d<3R$.
Its incoming and transverse derivative sets have fixed neighborhoods where $S=x\cdot d$ and $a_0=1$.
The phase and amplitude retain their eikonal and transport equations on a neighborhood of $\supp\chi_*$.
Put $b_0(\xi)=|\xi|^\mu-p_{\rm sm}(\xi)$ and set
\[
U=e^{iS/h}a_0,\qquad b_h=\chi_*W_h,\qquad u_{\rm app}=e^{ix\cdot d/h}+b_h.
\]
For $P_0=|hD|^\mu-1$, the exact identities are
\begin{equation}\label{eq:c4-exact-residual}
\begin{split}
E=P_{h,V}u_{\rm app}
&=P_0b_h+\delta Vu_{\rm app}\\
&=E_I+\mathcal C_h+b_0(hD)b_h,\\
E_I&=\chi_*(p_{\rm sm}(hD)-1+\delta V)U,\\
\mathcal C_h&=[p_{\rm sm}(hD),\chi_*]W_h.
\end{split}
\end{equation}
Apply the finite symbol expansion \eqref{eq:finite-symbolic-remainder} from Appendix~\ref{app:finite-kernel} with $L=2$.
The $C^2$ amplitude bound suffices, and the eikonal and transport equations cancel the first two coefficients, giving
\[
E_I=h^2\delta e^{iS/h}b_{I,h},\qquad \norm{b_{I,h}}_\infty\le C,
\qquad \supp b_{I,h}\subset\supp\chi_*.
\]

Choose a smooth transverse cutoff $\chi_T$ equal to one on a fixed neighborhood of the tube of radius $r_0$, and supported in the tube of radius $r_1$.
Write $E_{I,T}=\chi_T E_I$.
Then
\begin{equation}\label{eq:c4-interior-tube}
\norm{E_{I,T}}_2\le Ch^2\delta,
\qquad
\supp E_{I,T}\subset
\{y:|y-(y\cdot d)d|\le r_1,\ |y\cdot d|\le D\},
\end{equation}
where $D$ is fixed.
The remainder $(1-\chi_T)E_I$ has weighted $L^2$ and $L^\infty$ norms $O(h^M)$ for every fixed $M$.
On its support, $W_h$ and $V$ vanish locally, so the smooth-symbol kernel acts across a fixed positive transverse separation.
Its rapid decay away from the diagonal gives arbitrary powers of $h$ from the uniform $L^1$ bound of $W_h$.

Let $\chi_D\in C_c^\infty$ equal one near the downstream derivative support of $\chi_*$ and have support in $2R<x\cdot d<3R$.
Put $E_D=\chi_D\mathcal C_h$.
The incoming and transverse parts of $\mathcal C_h$ have arbitrary-order tails because $W_h$ vanishes on fixed neighborhoods there.
Hence we have the exact decomposition
\begin{equation}\label{eq:c4-source-decomposition}
E=E_{I,T}+E_D+E_T+b_0(hD)b_h,
\qquad
\norm{\langle x\rangle E_T}_2+\norm{E_T}_\infty\le C_Mh^M
\end{equation}
for every prescribed $M$.
The commutator expansion \eqref{eq:commutator-remainder} with $L=2$ gives $\norm{E_D}_1+\norm{E_D}_2\le Ch$ without additional spatial derivatives.

Fix $\overline{B_R}\subset\Omega_c\Subset B_{3R/2}$.
The downstream source is separated from this set in the $d$ direction.
Equation~\eqref{eq:backward-cone} gives, for every fixed $M$,
\begin{equation}\label{eq:c4-downstream-filtered}
\norm{R_{h,0}Q_hE_D}_{L^2(\Omega_c)}+\norm{R_{h,0}Q_hE_D}_{L^\infty(\Gamma)}\le C_Mh^M.
\end{equation}
Since the source $E_{I,T}$ satisfies the tube condition of Lemma~\ref{lem:directed-free-trace}, \eqref{eq:microlocal-sphere} gives
\begin{equation}\label{eq:c4-directional-trace-input}
\norm{R_{h,0}Q_hE_{I,T}}_{L^2(\Gamma)}\le Ch^{-1}\norm{E_{I,T}}_2\le Ch^{1+\mu}.
\end{equation}

\medskip
\noindent\textit{The complementary frequencies.}
For a compactly supported $C^2$ amplitude $a$, the finite expansion of $Q_h^\perp(e^{iS/h}a)$ has vanishing coefficients of orders zero and one.
Its kernel is a delta distribution minus a Schwartz kernel, and the remainder is bounded by
\begin{equation}\label{eq:c4-offphase-filter}
\norm{Q_h^\perp(e^{iS/h}a)}_\infty
+\norm{\langle x\rangle Q_h^\perp(e^{iS/h}a)}_2
\le Ch^2\norm a_{C^2}.
\end{equation}
Outside a fixed neighborhood of the amplitude support, it has arbitrary-order Schwartz tails.

Put $z_h=Q_h^\perp b_h$ and $F_h=Q_h^\perp(Vu_{\rm app})$.
For $\mu\le1$, apply \eqref{eq:c4-offphase-filter} separately to the two phases in $b_h$.
For $\mu\ge1$, write
\[
b_h=e^{ix\cdot d/h}\chi_*A_h,
\qquad A_h=e^{i(S-x\cdot d)/h}a_0-1,
\qquad \norm{A_h}_{C^2}\le C\delta/h.
\]
Since $Vu_{\rm app}=Ve^{iS/h}a_0$, the second filtered input has a compactly supported amplitude with a uniform $C^2$ bound.
Thus, with $\kappa_\mu=\max(2,1+\mu)$,
\begin{equation}\label{eq:c4-offphase-bounds}
\norm{z_h}_\infty+\norm{\langle x\rangle z_h}_2\le Ch^{\kappa_\mu},
\qquad
\norm{F_h}_\infty+\norm{\langle x\rangle F_h}_2\le Ch^2.
\end{equation}
For any fixed $r<r'<R$, the part of $F_h$ outside $B_{r'}$ has arbitrary-order weighted $L^2$ and pointwise bounds, by the Schwartz kernel of $q_d(hD)$ and the support of $Vu_{\rm app}$ in $B_r$.

The exact resolvent identities give
\begin{equation}\label{eq:c4-offphase-exact-identity}
\begin{split}
R_{h,0}Q_h^\perp E&=z_h+\delta R_{h,0}F_h,\\
-R_{h,V}Q_h^\perp E
&=-z_h+\delta R_{h,V}(Vz_h)-\delta R_{h,V}F_h.
\end{split}
\end{equation}
The second identity uses $R_{h,V}P_0z_h=z_h-\delta R_{h,V}(Vz_h)$.
Here $z_h\in H^2$ and is rapidly decreasing outside a fixed compact set.
Since $0<\mu<2$, the function $P_0z_h$ lies in weighted $L^2$; its fractional tail is $O(|x|^{-3-\mu})$.
The outgoing identities follow for smooth potentials and then by approximation.
The weighted resolvent estimate gives
\[
\norm{R_{h,V}Q_h^\perp E}_{L^2(B_R)}
\le C\bigl(h^{\kappa_\mu}+h^{\kappa_\mu+\mu-1}+h^{1+\mu}\bigr)
\le Ch^{1+\mu}.
\]
For the free sphere estimate, split $F_h$ with a smooth cutoff supported in $B_{r'}$ and equal to one near $B_r$.
The separated receiver estimate \eqref{eq:separated-sphere-L2} controls the compact part by $Ch^{-1}h^2$.
The weighted classical kernel bound and the fractional correction estimate \eqref{eq:surface-j-mass} control the remaining arbitrary-order tail.
Consequently
\begin{equation}\label{eq:c4-offphase-free-sphere}
\norm{R_{h,0}Q_h^\perp E}_{L^2(\Gamma)}
\le C\bigl(h^{\kappa_\mu}+h^{1+\mu}\bigr).
\end{equation}

\medskip
\noindent\textit{Comparison with the outgoing field.}
For the near frequency part, $Q_hb_0(hD)=0$ exactly.
Apply the weighted perturbed resolvent to $Q_hE_{I,T}$ and the arbitrary-order tail.
The former has weighted $L^2$ norm $Ch^2\delta$, since $Q_h$ is uniformly bounded on that space.
For the downstream term, use
\[
R_{h,V}Q_hE_D
=R_{h,0}Q_hE_D-\delta R_{h,V}\bigl(VR_{h,0}Q_hE_D\bigr)
\]
and \eqref{eq:c4-downstream-filtered}.
The error is $u-u_{\rm app}=-R_{h,V}E$.
These estimates, together with \eqref{eq:c4-offphase-exact-identity}, give
\begin{equation}\label{eq:c4-volume-final}
\norm{u-u_{\rm app}}_{L^2(B_R)}\le Ch^{1+\mu}.
\end{equation}
Since $u_{\rm app}=e^{iS/h}a_0$ near $\overline{B_R}$, this proves \eqref{eq:wkb-l2}.

For $V\in C_c^4(B_r)$, take a smooth approximation in $C^4$ supported in $B_r$.
The flow, momentum, and incoming label converge jointly in $C^3$, and the leading amplitude converges in $C^2$.
The two path identities used in \eqref{eq:phase-third-mixed} give convergence of $\partial_\delta^3S$ in $C^1_{x,d}$.
The finite kernel remainders converge in the bounded norms used above.
The localized Fredholm inverse is continuous at the common admissible high frequencies.
The constants depend only on the stated $C^4$ bound and fixed geometry, so the estimates pass to the limit.
\end{proof}

\subsection{Approximation on the receiver sphere}\label{sec:surface-approximation}

The volume error also controls the feedback from the potential support to the receiver sphere.

\begin{lemma}\label{lem:surface-wkb}
Under the hypotheses of Lemma~\ref{lem:wkb-l2}, the same phase and positive amplitude satisfy
\begin{equation}\label{eq:surface-wkb}
\sup_{d\in\Sph} \norm{u_{V,h,d}-e^{iS_{\delta,V,d}/h}a_{0,\delta,V,d}} _{L^2(\Gamma)} \le Ch^{\nu_\mu},\qquad \nu_\mu=\min(2\mu,1+\mu),\quad 0<\mu<2.
\end{equation}
The constants depend only on $\mu,r,R,F$.
\end{lemma}

\begin{proof}
Use the approximate field $u_{\rm app}$ and residual $E$ from the proof of Lemma~\ref{lem:wkb-l2}.
The integral equation gives
\[
u-u_{\rm app}=-R_{h,0}E-\delta R_{h,0}\bigl(V(u-u_{\rm app})\bigr).
\]
Use \eqref{eq:c4-directional-trace-input} on $E_{I,T}$, \eqref{eq:c4-downstream-filtered} on $E_D$, and \eqref{eq:c4-offphase-free-sphere} on the complementary part.
The Schwartz kernel of $Q_h$ preserves the weighted $L^2$ and $L^\infty$ tail bounds uniformly in $h$.
For the tails, Cauchy--Schwarz against the classical kernel and
\[
\sup_{x\in\Gamma}\int_{\R^3}\frac{\dd y}{\langle y\rangle^2|x-y|^2}<\infty
\]
give a factor $Ch^{-2}$ from the weighted $L^2$ norm.
The fractional correction is bounded by \eqref{eq:surface-j-mass} and the $L^\infty$ norm.
Choosing $M\ge5$ in \eqref{eq:c4-source-decomposition} bounds these terms by $Ch^3$.
The separated receiver estimate \eqref{eq:separated-sphere-L2} and the volume error \eqref{eq:c4-volume-final} control the feedback term by $Ch^{2\mu}$.
It follows that
\[
\norm{u-e^{iS/h}a_0}_{L^2(\Gamma)}
\le C\bigl(h^{1+\mu}+h^{\kappa_\mu}+h^{2\mu}+h^3\bigr)
\le Ch^{\min(2\mu,1+\mu)}.\qedhere
\]
\end{proof}

For $\mu\le1$, the squared spherical error is $O(h^{4\mu})$, the order used in the real-frequency estimate below.

\section{Increasing stability for the inverse scattering problem}\label{sec:stability}

In this section, we prove Theorem \ref{thm:main-stability} which establishes the increasing stability for the inverse potential scattering problem.
Let $W=V_1-V_2$ and retain $\mu=2\alpha$.
We use ray transforms and integer phase combinations to obtain a real-frequency estimate, then apply analytic continuation to choose the frequency window.

\subsection{Ray estimates and phase differences}\label{sec:ray-phase}

\begin{lemma}\label{lem:xray}
For every $f\in L^2(\R^3)$ supported in $B_r$ and either observation set in \eqref{eq:geometry},
\begin{equation}\label{eq:xray-norm}
\norm f_{H^{-1/2}(\R^3)}\le C\norm{Bf}_Y.
\end{equation}
\end{lemma}

\begin{proof}
Parameterize oriented lines by $(p,d)$, where $p\in d^\perp$.
The Fourier slice identity and Plancherel's theorem give
\begin{equation}\label{eq:slice}
\int_{\Sph}\int_{d^\perp}|Xf(p,d)|^2\dd p\dd d =c\int_{\R^3}|\xi|^{-1}|\widehat f(\xi)|^2\dd\xi \ge c'\norm f_{H^{-1/2}}^2.
\end{equation}
Compact support and $L^2$ regularity give $f\in L^1$, so the integral is finite near the origin.

For spherical observations, parameterize the downstream hemisphere by $x=p+\sqrt{R^2-|p|^2}\,d$.
The half-ray integral equals $Xf(p,d)$ at this point and vanishes at the upstream intersection.
Therefore,
\begin{equation}\label{eq:surface-xray}
\int_\Gamma|Bf(x,d)|^2\dd\sigma(x) =\int_{|p|<r}\frac R{\sqrt{R^2-|p|^2}} |Xf(p,d)|^2\dd p.
\end{equation}
The Jacobian is bounded above and below on $|p|<r<R$.
Combining this identity with \eqref{eq:slice} proves \eqref{eq:xray-norm} for $\mathcal O=\Gamma$.

For annular observations, the part of a line in $A$ has length
\[
\ell(p)=2\bigl(\sqrt{R^2-|p|^2}-\sqrt{r^2-|p|^2}\bigr), \qquad 2(R-r)\le\ell(p)\le2R
\]
when $|p|<r$.
Since $Xf$ is constant on each line and vanishes for $|p|\ge r$, Fubini's theorem gives
\[
\norm{Xf}_{L^2(A\times\Sph)}^2 =\int_{\Sph}\int_{d^\perp}\ell(p)|Xf(p,d)|^2\dd p\dd d.
\]
The result follows from \eqref{eq:slice} and $Xf=Bf(\cdot,d)+Bf(\cdot,-d)$.
\end{proof}

We next compare the phases associated with the two potentials.
Set $G(k)=u_1^s(k)-u_2^s(k)$, $\beta=BW$, and
\[
\phi_j(x,d,k)=k(S_{k^{-\mu},V_j,d}(x)-x\cdot d), \qquad \psi=\phi_1-\phi_2.
\]
The joint regularity in the independent parameter $\delta$ yields the following phase expansion.

\begin{lemma}\label{lem:phase-taylor}
There is a real function $c_2(x,d)$, independent of $k$, such that, for sufficiently large $k$,
\begin{equation}\label{eq:phase-taylor}
\begin{split}
\psi(k)&=-\frac{\beta}{\mu}k^{1-\mu} +c_2k^{1-2\mu}+\mathfrak r(k),\\
\norm{\mathfrak r(k)}_{C^1}&\le Ck^{1-3\mu},\qquad \norm{\partial_k\mathfrak r(k)}_{C^1}\le Ck^{-3\mu}.
\end{split}
\end{equation}
Here the $C^1$ norm is taken in the spatial and incident-direction variables on the fixed construction sets.
\end{lemma}

\begin{proof}
The mixed bound \eqref{eq:phase-third-mixed} in the phase construction gives a uniformly bounded third $\delta$ derivative of the phase in $C^1_{x,d}$.
For $D_\delta=S_{\delta,V_1,d}-S_{\delta,V_2,d}$, Taylor's formula therefore gives
\[
D_\delta=-\delta\beta/\mu+\delta^2c_2+r_\delta, \qquad c_2=\tfrac12\partial_\delta^2D_0, \qquad r_\delta=\tfrac12\int_0^\delta (\delta-t)^2\partial_\delta^3D_t\dd t.
\]
Thus $\norm{r_\delta}_{C^1}\le C\delta^3$ and $\norm{\partial_\delta r_\delta}_{C^1}\le C\delta^2$.
Taking $\mathfrak r(k)=kr_{k^{-\mu}}$ and using $\partial_k\mathfrak r=r_\delta-\mu\delta\partial_\delta r_\delta$ proves the assertion.
\end{proof}

Positive leading amplitudes yield the following estimate for integer combinations of phases.

\begin{lemma}\label{lem:integer-phase}
Let $t_1,\ldots,t_m>0$ and $n_1,\ldots,n_m\in\mathbb Z$ be fixed, and set $\Psi(k)=\sum_{\ell=1}^m n_\ell\psi(t_\ell k)$.
Put $\nu_\mu=\min(2\mu,1+\mu)$.
Then, uniformly in $d$ and for sufficiently large $k$,
\begin{equation}\label{eq:integer-phase-data}
\norm{e^{i\Psi(k)}-1}_{L^2(\mathcal O)}^2 \le C\left(\sum_{\ell=1}^m \norm{G(\cdot,d,t_\ell k)}_{L^2(\mathcal O)}^2 +k^{-2\nu_\mu}\right).
\end{equation}
The same inequality holds in $Y$ after integration in $d$.
The constant and the frequency threshold may also depend on the fixed ratios $t_\ell$ and integers $n_\ell$.
\end{lemma}

\begin{proof}
Lemmas~\ref{lem:wkb-l2} and \ref{lem:surface-wkb} give $e^{-ikx\cdot d}u_j=e^{i\phi_j}(a_j+\rho_j)$, with $a_j\ge3/4$ real and $\norm{\rho_j(\cdot,d,k)}_{L^2(\mathcal O)}\le Ck^{-\nu_\mu}$.
The identity
\[
|a_1e^{i\phi_1}-a_2e^{i\phi_2}|^2 =(a_1-a_2)^2+a_1a_2|e^{i\psi}-1|^2
\]
implies $|e^{i\psi(k)}-1|^2 \le C(|G(x,d,k)|^2+|\rho_1|^2+|\rho_2|^2)$ almost everywhere.
Telescoping on the unit circle gives
\[
|e^{i\Psi(k)}-1| \le\sum_{\ell=1}^m|n_\ell| |e^{i\psi(t_\ell k)}-1|.
\]
Squaring and integrating proves \eqref{eq:integer-phase-data}.
\end{proof}

\subsection{The real-frequency estimate}\label{sec:real-frequency}

Integer combinations of phases remove the second coefficient in \eqref{eq:phase-taylor}.
For $\mu<1$, frequency averaging controls the resulting oscillations; at $\mu=1$, the incoming cap fixes the phase branch; for $\mu>1$, the phase tends to zero.
These arguments give a uniform estimate on a fixed relative window.

\begin{proposition}\label{prop:filtered-real-window}
Set $\tau=(\mu-1)_+$ and $\rho=\min(4\mu,4)$.
There are $K_*,C>0$, depending only on $\alpha,r,R,F$, such that, for every $Q\ge4K_*$,
\begin{equation}\label{eq:filtered-real-window}
\norm W_{H^{-1/2}(\R^3)}^2 \le C\left(\frac1Q\int_{Q/4}^Q t^{2\tau}\norm{G(t)}_Y^2\dd t+Q^{-\rho}\right).
\end{equation}
\end{proposition}

\begin{proof}
By Lemma~\ref{lem:xray}, it suffices to estimate $\norm\beta_Y^2$.
We treat the three ranges of $\mu$ separately.

\medskip
\noindent\textit{The case $0<\mu<1$.}
Take $a=3/4$, $p=1-2\mu$, and
\[
\lambda=\begin{cases}
(2-a^p)^{1/p},&p\ne0,\\
4/3,&p=0.
\end{cases}
\]
The power mean of order $p$ of $a,\lambda$ equals one.
Comparison with the orders $-1$ and $1$ gives $5/4\le\lambda\le3/2$.
Since $1-\mu>p$, strict monotonicity of power means gives $c_\mu=a^{1-\mu}+\lambda^{1-\mu}-2>0$.
Thus $\Psi(k)=\psi(ak)-2\psi(k)+\psi(\lambda k)$ cancels the second coefficient in \eqref{eq:phase-taylor}, and
\begin{equation}\label{eq:subcritical-filter-derivative}
\Psi'(k)=-\frac{c_\mu(1-\mu)}{\mu} \beta k^{-\mu}+O(k^{-3\mu}).
\end{equation}

Fix $(x,d)$ and put $I_Q=[Q/(2\lambda),Q/\lambda]$.
We have $|\beta(x,d)|\le M_B$ uniformly.
If $|\beta|\le A_0Q^{-2\mu}$, its square is bounded by $A_0^2Q^{-4\mu}$.
Otherwise choose the fixed $A_0$ large enough that \eqref{eq:subcritical-filter-derivative} yields
\[
c|\beta|Q^{-\mu}\le|\Psi'(k)|\le C|\beta|Q^{-\mu} \quad(k\in I_Q),
\]
with a fixed sign.
The phase is therefore monotone on $I_Q$, with image length comparable to $T=Q^{1-\mu}|\beta|$.
For every real interval $J$ of length $\ell$,
\[
\frac1\ell\int_J|e^{is}-1|^2\dd s \ge c\min(1,\ell^2);
\]
indeed, the minimum over its location is $2-4|\sin(\ell/2)|/\ell$.
Changing variables to $s=\Psi(k)$ gives
\[
\frac1Q\int_{I_Q}|e^{i\Psi(k)}-1|^2\dd k \ge c\min(1,T^2)\ge c_{M_B}|\beta|^2.
\]
Including the small-$\beta$ case, integrating in $(x,d)$, and applying Lemma~\ref{lem:integer-phase}, we obtain
\[
\norm\beta_Y^2\le C\left( \frac1Q\int_{Q/4}^Q\norm{G(t)}_Y^2\dd t +Q^{-4\mu}\right).
\]
Here $aI_Q,I_Q,\lambda I_Q\subset[Q/4,Q]$, since $a/(2\lambda)\ge1/4$.
Since $2\nu_\mu=4\mu$, the field approximation contributes the same remainder order.

\medskip
\noindent\textit{The case $\mu=1$.}
Set $\Psi(k)=2\psi(2k)-\psi(k)$.
Lemma~\ref{lem:phase-taylor} gives
\begin{equation}\label{eq:critical-filter-phase}
\Psi(k)=-\beta+O_{C^1}(k^{-2}).
\end{equation}
Fix $d$, put $Z=e^{i\Psi(k)}-1$, and write
\[
e_d^2=\norm{G(\cdot,d,k)}_{L^2(\mathcal O)}^2 +\norm{G(\cdot,d,2k)}_{L^2(\mathcal O)}^2.
\]
The spatial Lipschitz norm of $Z$ is uniformly bounded, and Lemma~\ref{lem:integer-phase} gives $\norm Z_{L^2(\mathcal O)}\le C(e_d+k^{-2})$.
The observation measure satisfies $m_{\mathcal O}(\mathcal O\cap B(x,t))\ge ct^{n_{\mathcal O}}$ for small $t$, where $n_{\mathcal O}=2$ on the sphere and $3$ on the annulus.
Integrating on a ball whose radius is proportional to the maximum of $|Z|$ therefore yields
\[
\norm Z_\infty\le C(e_d+k^{-2})^{2/(n_{\mathcal O}+2)}.
\]
Choose a fixed small $e_*>0$ and then increase $K_*$.
If $e_d\le e_*$, the last bound makes $\norm Z_\infty<1$.
Since the observation set is connected and $\Psi=0$ on its nonempty incoming cap, the phase lies in the component containing zero of $\{s:|e^{is}-1|<1\}$.
Hence $|\Psi|\le C|Z|$ on that component.
Equation~\eqref{eq:critical-filter-phase} gives
\[
\norm{\beta(\cdot,d)}_{L^2(\mathcal O)}^2 \le C(e_d^2+k^{-4}).
\]
For the remaining directions, $e_d>e_*$, the uniform a priori bound for $\beta$ gives the same inequality.
Integrating in $d$ and then averaging over $k\in[Q/4,Q/2]$ proves
\[
\norm\beta_Y^2\le C\left( \frac1Q\int_{Q/4}^Q\norm{G(t)}_Y^2\dd t+Q^{-4}\right).
\]

\medskip
\noindent\textit{The case $1<\mu<2$.}
Take $\lambda=2^{1/(2\mu-1)}\in(1,2)$ and $\Psi(k)=2\psi(\lambda k)-\psi(k)$.
Since $2\lambda^{1-2\mu}=1$, we have
\[
\begin{split}
\Psi(k)&=-\frac{d_\mu\beta}{\mu}k^{1-\mu} +O(k^{1-3\mu}),\\
d_\mu&=2\lambda^{1-\mu}-1 =2^{\mu/(2\mu-1)}-1>0.
\end{split}
\]
Since $\norm\Psi_\infty\le Ck^{1-\mu}\to0$, we have $|\Psi|\le C|e^{i\Psi}-1|$ for sufficiently large $k$.
Lemma~\ref{lem:integer-phase} consequently gives
\[
\norm\beta_Y \le C\left(k^{\mu-1} \bigl(\norm{G(k)}_Y+\norm{G(\lambda k)}_Y\bigr) +k^{-2}+k^{-2\mu}\right).
\]
The last term is bounded by $k^{-2}$.
Squaring and averaging over $k\in[Q/(2\lambda),Q/\lambda]$, for which $k,\lambda k\in[Q/4,Q]$, gives
\[
\norm\beta_Y^2\le C\left( \frac1Q\int_{Q/4}^Q t^{2(\mu-1)} \norm{G(t)}_Y^2\dd t+Q^{-4}\right).
\]
Finally, Lemma~\ref{lem:xray} proves \eqref{eq:filtered-real-window} in all three cases.
Choose $K_*$ large enough that every sampled frequency satisfies the geometric optics and phase thresholds.
\end{proof}

\subsection{Analytic continuation of the measured fields}\label{sec:extended-domain}

We extend the measured fields from the fixed strip to a domain where quantitative continuation controls the data beyond the measured interval.

\begin{proposition}\label{prop:extended-domain}
There are fixed $b>0$ and $L>b+1$, uniform for $V\in\mathcal V_{r,F}^4$, such that the localized outgoing resolvent and the $Y$-valued scattered field extend holomorphically to
\begin{equation}\label{eq:extended-domain}
\mathcal D_{b,L} =\{z\in\C:\Im z>-b,\ |z+ib|>L\}.
\end{equation}
Use the branch of $\log z$ with cut along $-i[0,\infty)$ and argument in $(-\pi/2,3\pi/2)$.
For fixed compactly supported smooth cutoffs $\eta_1,\eta_2$,
\begin{equation}\label{eq:extended-resolvent-bound}
\norm{\eta_1\cR_{\alpha,V}(z)\eta_2}_{L^2\to L^2} \le C|z|^{1-\mu},\qquad z\in\mathcal D_{b,L}.
\end{equation}
Moreover,
\begin{equation}\label{eq:extended-data-bound}
\norm{u_V^s(z)}_Y \le C(1+|z|)^3 e^{r(\Im z)_+}.
\end{equation}
These extensions agree with the outgoing fields on the sufficiently large positive real axis.
\end{proposition}

The spectral and outgoing resolvents agree on their first-quadrant overlap, and a spectral Sobolev estimate provides the spherical trace.
The proof is given in Appendix~\ref{app:extended-domain}.

\begin{corollary}\label{cor:normalized-discrepancy}
Let $G(z)=u_{V_1}^s(z)-u_{V_2}^s(z)$ and $\tau=(\mu-1)_+$.
Fix $\kappa>b+1$ and an integer $N\ge4$.
Then
\begin{equation}\label{eq:normalized-discrepancy}
F_N(z)=e^{irz}(z+i\kappa)^{-N}z^\tau G(z)
\end{equation}
is holomorphic and uniformly bounded in $Y$ on $\mathcal D_{b,L}$.
Its norm on a measured positive real interval is at most the weighted discrepancy there, and
\begin{equation}\label{eq:undo-normalization}
\norm{q^\tau G(q)}_Y \le C(1+q)^N\norm{F_N(q)}_Y
\end{equation}
for positive real $q$ in the domain.
\end{corollary}

\begin{proof}
The branch defining $z^\tau$ is holomorphic in the domain, and the pole $-i\kappa$ lies below it.
Moreover,
\[
|z+i\kappa|\ge\kappa-b>1,\qquad 1+|z|\le C_{\kappa,b}|z+i\kappa|.
\]
The factor $|e^{irz}|=e^{-r\Im z}$ cancels the exponential growth in \eqref{eq:extended-data-bound} and is bounded below the real axis.
Since $\tau<1$ and $N\ge4$, this proves the uniform bound.
On the positive real axis the exponential has modulus one and $|q+i\kappa|^{-N}\le1$, which proves the remaining assertions.
\end{proof}

\subsection{Choice of the frequency window}

In what follows, we apply the analytic continuation given by Lemma~\ref{lem:exterior-continuation} to the normalized field and balance the continued data against the real-frequency remainder,
which yields the proof of Theorem~\ref{thm:main-stability}.

\begin{proof}[Proof of Theorem~\ref{thm:main-stability}]
Put $\mu=2\alpha$, $\rho=\min(4\mu,4)$, $G(z)=u_1^s(z)-u_2^s(z)$, and $H(z)=z^\tau G(z)$.
Fix $N=4$ in Corollary~\ref{cor:normalized-discrepancy} and choose $M\ge1$, depending only on the a priori parameters, so that
\[
F(z)=M^{-1}e^{irz}(z+i\kappa)^{-N}H(z)
\]
is bounded by one in $Y$ on $\mathcal D_{b,L}$.
Choose $K_0\ge2L+2$ sufficiently large for all real-frequency estimates.
On the measured interval, $\norm{F(k)}_Y\le\eps_w$.
Lemma~\ref{lem:exterior-continuation} and \eqref{eq:undo-normalization} therefore yield
\begin{equation}\label{eq:continued-weighted-data}
\norm{H(q)}_Y\le C(1+q)^N \exp\left(-\frac{cE}{(1+q)^2}\right), \quad q\ge K_0+1,\qquad E=|\log\eps_w|.
\end{equation}

Set
\[
\mathcal L(E)=\sqrt{\frac{E}{\log(2+E)}}.
\]
Choose $0<\sigma\le1$ so small that
\begin{equation}\label{eq:continuation-sigma}
\frac{c}{4\sigma^2}> \frac N2+\frac\rho4+1, \qquad a=\frac\sigma8.
\end{equation}
If $K\ge a\mathcal L(E)$ and $K\ge4K_0$, apply Proposition~\ref{prop:filtered-real-window} with $Q=K$.
Every frequency in $[K/4,K]$ is measured, and hence
\[
\norm W_{H^{-1/2}}^2 \le C(\eps_w^2+K^{-\rho}) \le C\left(\eps_w^2+ \bigl(K+\mathcal L(E)\bigr)^{-\rho}\right).
\]
If $K\ge a\mathcal L(E)$ but $K<4K_0$, then $K+\mathcal L(E)\le4K_0(1+a^{-1})$.
The a priori bound gives the same estimate.

Suppose finally that $K<a\mathcal L(E)$ and take $Q=\sigma\mathcal L(E)$.
Then
\[
Q/4=2a\mathcal L(E)>2K\ge2(K_0+1),
\]
so $[Q/4,Q]$ lies beyond the measured band and above the fixed frequency thresholds.
For each $k$ in this interval, \eqref{eq:continued-weighted-data} gives
\[
\norm{H(k)}_Y \le C Q^N\exp\left(-\frac{cE}{4Q^2}\right) =C Q^N(2+E)^{-c/(4\sigma^2)} \le C Q^{-\rho/2}.
\]
The last inequality follows from $Q\le E^{1/2}$ and \eqref{eq:continuation-sigma}.
Substituting this bound into Proposition~\ref{prop:filtered-real-window} yields
\[
\norm W_{H^{-1/2}}^2\le C Q^{-\rho} \le C\bigl(K+\mathcal L(E)\bigr)^{-\rho},
\]
because $Q\asymp K+\mathcal L(E)$ in the present case.
This proves \eqref{eq:main-stability}.

For $0<\theta<1/2$, the ratio $\mathcal L(E)/E^\theta$ has a positive infimum on $[1,\infty)$ and tends to infinity as $E\to\infty$.
Consequently $K+\mathcal L(E)\ge c_\theta(K+E^\theta)$.
The a priori $H^s$ bound for $0<s\le4$ and interpolation give
\[
\norm W_{L^2}^2 \le \norm W_{H^{-1/2}}^{4s/(2s+1)} \norm W_{H^s}^{2/(2s+1)}.
\]
This proves the $L^2$ conclusion. The proof is complete.
\end{proof}

% Numerical section based exclusively on the audited 72-receiver experiment.
% Required packages: graphicx, booktabs.
\section{Numerical experiments}\label{sec:dense-numerics}

In this section, we present numerical examples to examine the effect of
the measured frequency band on the reconstruction.
For ease of visualization, we carry out the experiments in two dimensions.

We consider the two-dimensional fractional Helmholtz equation,
\begin{equation}\label{eq:dense-model}
 (-\Delta)^\alpha u-k^{2\alpha}u+Vu=0
 \quad\text{in }\mathbb R^2,\qquad
 u=e^{ikx\cdot d}+u^s,\quad d\in\mathbb S^1.
\end{equation}
The complex scattered field is measured on a circle enclosing the
support of $V$. The examples cover the three regimes
$\alpha<1/2$, $\alpha=1/2$, and $\alpha>1/2$.

\subsection{Reconstruction method}

Let $D=(-1,1)^2$ and let $G_{\alpha,k}^{(2)}$ be the outgoing Green
function of $(-\Delta)^\alpha-k^{2\alpha}$ in $\mathbb R^2$.
Writing $u(x,d,k)=e^{ikx\cdot d}w(x,d,k)$ gives
\begin{equation}\label{eq:dense-state}
 w(x)+\int_D G_{\alpha,k}^{(2)}(x-y)e^{-ikd\cdot(x-y)}
                 V(y)w(y)\,\mathrm dy=1.
\end{equation}
Factoring out the incident plane wave from the total field removes the known incident phase while retaining the full scattering interaction. On a square grid, let $T_{j\ell}$ and
$H_{j\ell}$ be the cell integrals of the kernel in
\eqref{eq:dense-state}, evaluated at interior points and receivers,
respectively. For $c\in\mathbb R^{81}$ and
$V_c=\sum_{a=1}^{81}c_a\psi_a$, set
$v=\Phi c$, where $\Phi_{na}=\psi_a(x_n)$. The discrete forward map is
\begin{equation}\label{eq:dense-forward}
 A_{j\ell}=I+T_{j\ell}\operatorname{diag}(v),\qquad
 A_{j\ell}w_{j\ell}=\boldsymbol1,\qquad
 f_{j\ell}(c)=-H_{j\ell}(v\odot w_{j\ell}).
\end{equation}
Here $\odot$ denotes entrywise multiplication. The entries of
$f_{j\ell}$ are the scattered field multiplied by
$e^{-ik_jz_m\cdot d_\ell}$; the measured data undergo the same
known phase change.

Let $\widetilde F_j(c)$ collect these columns and let
$\widetilde Y_j^\eta$ denote the noisy data. At upper endpoint $K$,
we minimize
\begin{equation}\label{eq:dense-objective}
 J_K(c)=
 \frac{\sum_{k_j\le K}\omega_{j,K}
       \|\widetilde F_j(c)-\widetilde Y_j^\eta\|_F^2}
      {2(K-2)N_rN_d\,s_\alpha^2}
       +\frac{\lambda_K}{2}|c|^2,\qquad
 \lambda_K=\frac{3.6\times10^{-4}}{K-2}.
\end{equation}
The subscript $F$ denotes the Frobenius norm,
$N_r=72$, $N_d=36$, and $\omega_{j,K}$ are trapezoidal weights
with sum $K-2$. The scale $s_\alpha$ is the root-mean-square
value of the noisy $k=2$ data for $V_1$, fixed for all examples
at that order. The penalty in the unnormalized cumulative functional
is therefore independent of $K$.

For comparison, we also compute regularized Born reconstructions by replacing the total field in the scattering integral with the incident field.
We use the same data, basis, weights, and regularization parameter to assess the effect of this linearization on reconstruction accuracy.

The gradient is evaluated by an adjoint solve. For
$r_{j\ell}=f_{j\ell}-\widetilde Y_{j,\ell}^\eta$ and
$\gamma_{j,K}=\omega_{j,K}/((K-2)N_rN_ds_\alpha^2)$, set
$b_{j\ell}=\gamma_{j,K}H_{j\ell}^*r_{j\ell}$ and solve
$A_{j\ell}^*p_{j\ell}=\operatorname{diag}(v)b_{j\ell}$. Then
\begin{equation}\label{eq:dense-gradient}
 \nabla J_K(c)=\lambda_Kc-
 \sum_{k_j\le K}\sum_{\ell=1}^{N_d}
 \Phi^T\operatorname{Re}
 \bigl[\overline{w_{j\ell}}\odot
       (b_{j\ell}-T_{j\ell}^*p_{j\ell})\bigr].
\end{equation}
We use BFGS with Armijo backtracking in coordinates $\zeta=L_K^Tc$,
where $L_KL_K^T$ is the regularized Born normal matrix.
The states and gradients remain nonlinear.
The initial vector is zero at $K=5$; each subsequent
band starts from the preceding reconstruction. Iterations stop when
the scaled gradient is below $10^{-6}$, or when the objective
decrease is below $2\times10^{-11}$ together with a small step
or gradient. At most three passes of $160$ iterations are allowed.

\subsection{Experimental setting}

The receivers are
$z_m=2.2(\cos\theta_m,\sin\theta_m)$, where
$\theta_m=2\pi(m+1/2)/72$, $m=0,\ldots,71$.
The incident angles are $2\pi\ell/36$, $\ell=0,\ldots,35$.
We sample the wavenumbers
$2,3,4,5,6,10,15,20,30,40,50$ and use the cumulative endpoints
$K=5,10,20,30,40,50$.
At every sampled wavenumber, synthetic data are generated on
$192^2$ cells and inversion uses $128^2$ cells.
Cell integrals include the oscillatory factor in
\eqref{eq:dense-state}; singular self-interactions are integrated
by a Duffy transformation. The interior operator is applied by
zero-padded FFTs and the linear systems are solved by GMRES.
Its relative tolerance is $2\times10^{-8}$ for data generation
and inversion, except that inversion at $\alpha=0.25$ uses
$2\times10^{-10}$.

For the reconstruction space, we take
\[
 \psi_a(x)=\chi(x)\exp\!\left(-\frac{|x-\xi_a|^2}{2(0.16)^2}\right),
 \qquad \xi_a\in\{-0.64,-0.48,\ldots,0.64\}^2.
\]
The radial cutoff is $\chi(x)=s((0.94-|x|)/0.16)$, where
$s(t)=\rho(t)/(\rho(t)+\rho(1-t))$ and
$\rho(t)=e^{-1/t}$ for $t>0$, with $\rho(t)=0$ otherwise.
Thus $\chi=1$ on $B_{0.78}$ and vanishes outside $B_{0.94}$.
To define the test potentials, put
\[
 g(x;a,b,\vartheta)=
 \exp\!\left(-\frac12\sum_{i=1}^2
       \frac{[R_\vartheta^T(x-a)]_i^2}{b_i^2}\right),
\]
where $R_\vartheta$ is the rotation matrix. We choose
\begin{align*}
 V_1(x)&=\chi(x)g(x;(0.12,-0.08),(0.24,0.34),\pi/7),\\
 V_2(x)&=\chi(x)\bigl[g(x;(-0.35,-0.18),(0.22,0.22),0)\\
       &\hspace{49pt}+0.8g(x;(0.35,0.20),(0.24,0.24),0)\bigr],\\
 V_3(x)&=\chi(x)\bigl[1.1g(x;(-0.28,0.20),(0.28,0.14),\pi/6)\\
       &\hspace{49pt}-0.85g(x;(0.32,-0.20),(0.20,0.25),-\pi/5)\bigr].
\end{align*}
These are a positive peak, two positive peaks, and a sign-changing
pair. All three are smooth and compactly supported.

At each frequency, we perturb the clean data $Y_j$ by
\begin{equation}\label{eq:dense-noise}
 Y_j^\eta=Y_j+\eta\|Y_j\|_F\frac{E_j}{\|E_j\|_F},
 \qquad\eta=0.03,
\end{equation}
where the real and imaginary parts of $E_j$ are independent
standard Gaussian arrays. One fixed realization is used for each
order, potential and frequency, and is retained when the band
expands. All reconstruction errors are evaluated on a common
$161^2$ grid using
$e_V(K)=\|V_{c_K}-V\|_{L^2(D)}/\|V\|_{L^2(D)}$.

% Keep spacing adjustments local to the numerical figures and tables.
\begingroup
\setlength{\floatsep}{8pt plus 2pt minus 2pt}
\setlength{\textfloatsep}{10pt plus 2pt minus 2pt}
\setlength{\intextsep}{10pt plus 2pt minus 2pt}
\subsection{Reconstruction results}

Table~\ref{tab:dense-errors} gives all relative reconstruction errors.
Figures~\ref{fig:alpha025-paired}--\ref{fig:alpha075-paired} show the
reconstructions at all six endpoints, together with their absolute errors.
The final errors range from $0.2448\%$ to $1.2564\%$; the corresponding
data residuals are $2.997\%$--$3.046\%$.

\begin{table}[!p]
 \centering
 \caption{Relative $L^2$ reconstruction errors (\%) with $3\%$ noise.}
 \label{tab:dense-errors}
 \small
 \setlength{\tabcolsep}{7pt}
 \begin{tabular}{ccrrrrrr}
 \toprule
 $\alpha$ & Potential & $K=5$ & $10$ & $20$ & $30$ & $40$ & $50$\\
 \midrule
 $0.25$ & $V_1$ & 1.6153 & 0.6041 & 0.3128 & 0.2584 & 0.2504 & 0.2448\\
        & $V_2$ & 2.8857 & 0.4750 & 0.3116 & 0.2910 & 0.2829 & 0.2709\\
        & $V_3$ & 7.6540 & 1.0696 & 0.6988 & 0.6616 & 0.6543 & 0.6489\\
 \addlinespace
 $0.50$ & $V_1$ & 0.8385 & 0.4980 & 0.3729 & 0.3578 & 0.3679 & 0.3700\\
        & $V_2$ & 1.2799 & 0.5240 & 0.4365 & 0.4189 & 0.4023 & 0.3870\\
        & $V_3$ & 8.0690 & 1.2532 & 0.8587 & 0.8054 & 0.7724 & 0.7420\\
 \addlinespace
 $0.75$ & $V_1$ & 0.8639 & 0.5202 & 0.4087 & 0.3715 & 0.3608 & 0.3534\\
        & $V_2$ & 1.3511 & 0.4903 & 0.4336 & 0.3932 & 0.3790 & 0.3661\\
        & $V_3$ & 9.1701 & 2.6255 & 1.5397 & 1.3840 & 1.3054 & 1.2564\\
 \bottomrule
 \end{tabular}
\end{table}

\begin{figure}[!htbp]
 \centering
 \includegraphics[width=\linewidth]{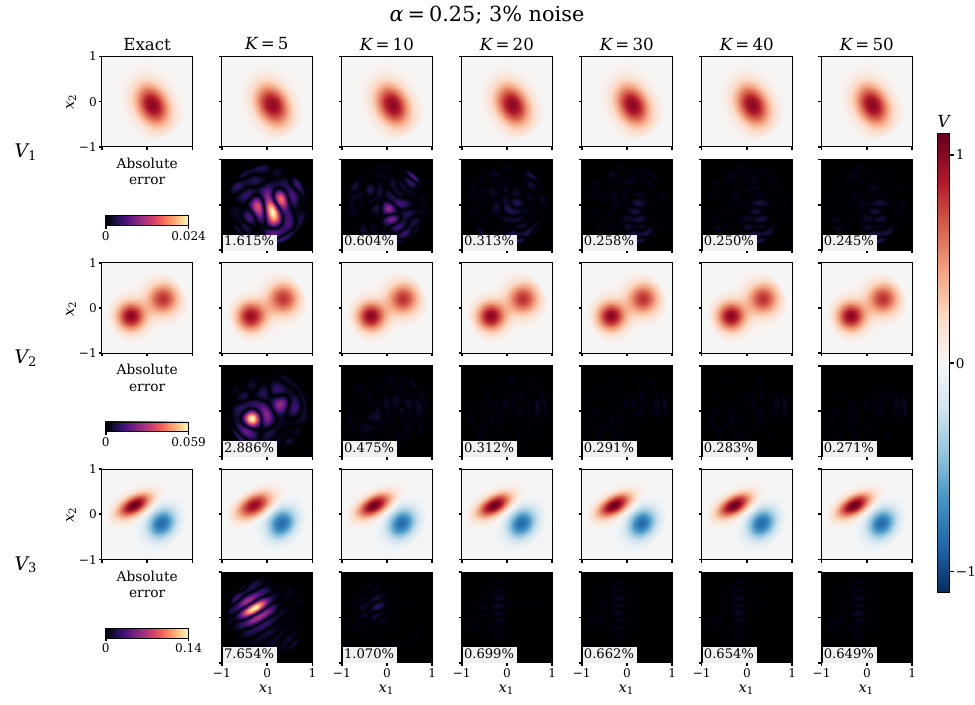}
 \caption{Reconstructions and absolute errors for $\alpha=0.25$.
 Row pairs correspond to $V_1,V_2,V_3$, with reconstruction above error.
 The first column contains exact potentials and error colorbars;
 the remaining columns show $K=5,10,20,30,40,50$.
 The right colorbar applies to potentials.
 For each potential, the error scale is fixed across all $K$ and $\alpha$.
 Labels give relative $L^2$ errors (\%).}
 \label{fig:alpha025-paired}
\end{figure}
\FloatBarrier

All nine examples improve between $K=5$ and $K=50$; eight curves
decrease monotonically. The exception is $V_1$ at $\alpha=0.5$,
whose error rises from $0.3578\%$ at $K=30$ to $0.3700\%$ at $K=50$.
The computed projection errors on the evaluation grid are
$0.1833\%$, $0.2053\%$, and $0.6198\%$ for $V_1,V_2,V_3$.
The final $V_3$ error $0.6489\%$ at $\alpha=0.25$ is close to this level.
Figure~\ref{fig:dense-bandwidth} collects these error curves and projection levels.

Table~\ref{tab:dense-born} compares nonlinear and Born reconstructions
of $V_3$ at $K=50$. The Born fits set $w_{j\ell}=\boldsymbol1$
in \eqref{eq:dense-forward}, with the same data, basis, weights and
penalty. 

\begin{table}[!htbp]
 \centering
 \caption{Born and nonlinear reconstruction errors (\%) for $V_3$ at $K=50$.}
 \label{tab:dense-born}
 \begin{tabular}{crr}
 \toprule
 $\alpha$ & Born & Nonlinear\\
 \midrule
 $0.25$ & 104.696 & 0.649\\
 $0.50$ &   6.472 & 0.742\\
 $0.75$ &   2.868 & 1.256\\
 \bottomrule
 \end{tabular}
\end{table}

\begin{figure}[!htbp]
 \centering
 \includegraphics[width=\linewidth]{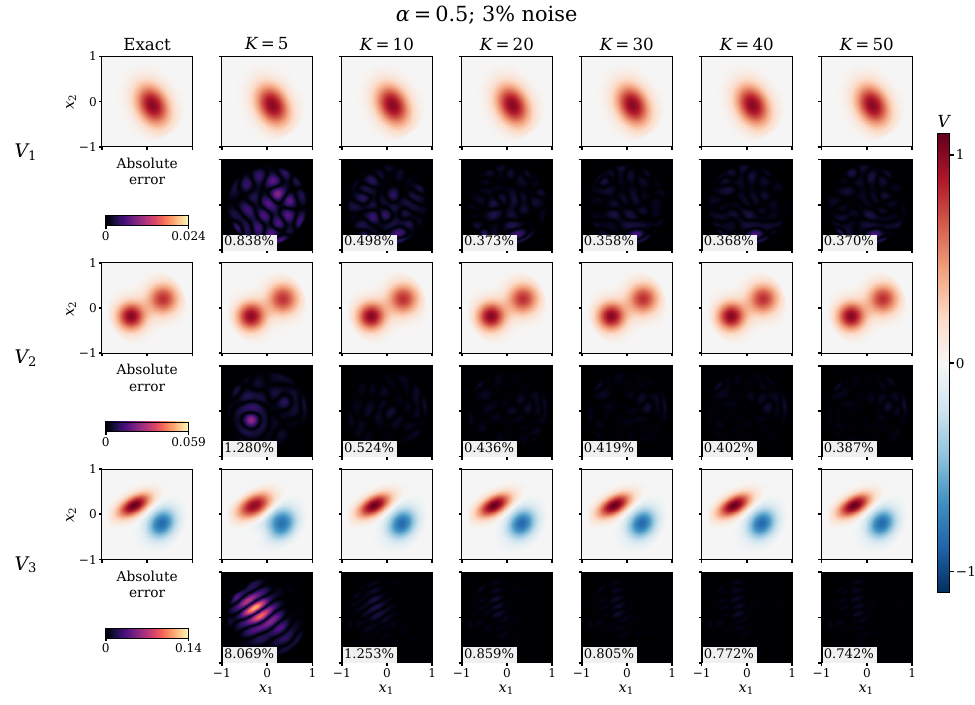}
 \caption{Reconstructions and absolute errors for $\alpha=0.5$,
 with the layout and color scales of Figure~\ref{fig:alpha025-paired}.}
 \label{fig:alpha050-paired}
\end{figure}

\clearpage
\begin{figure}[!htbp]
 \centering
 \includegraphics[width=\linewidth]{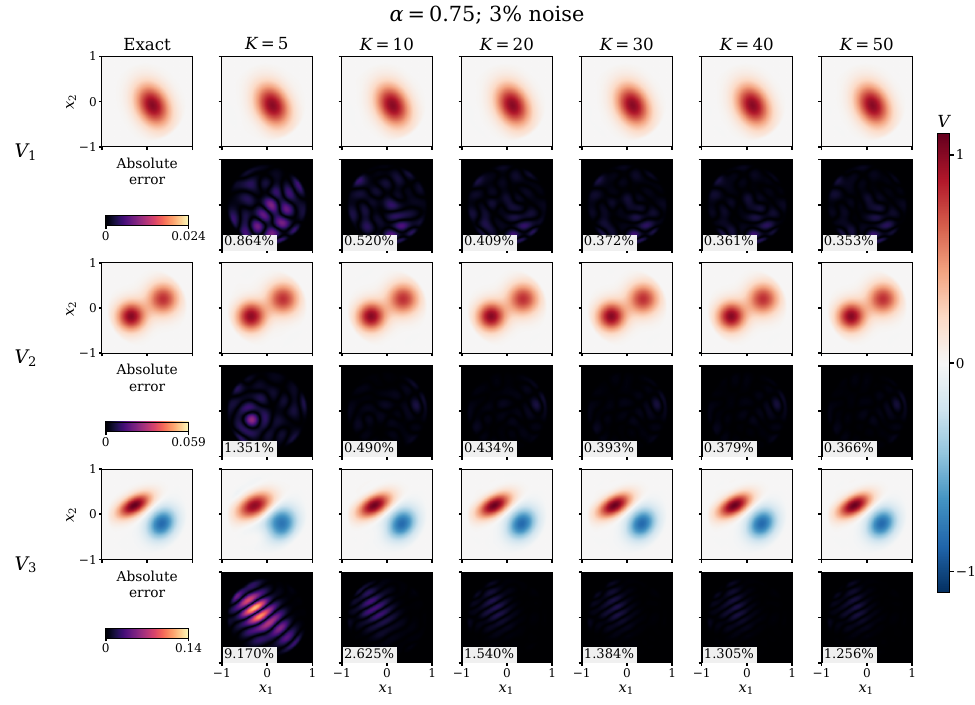}
 \caption{Reconstructions and absolute errors for $\alpha=0.75$,
 with the layout and color scales of Figure~\ref{fig:alpha025-paired}.}
 \label{fig:alpha075-paired}
\end{figure}

\begin{figure}[!htbp]
 \centering
 \includegraphics[width=\linewidth]{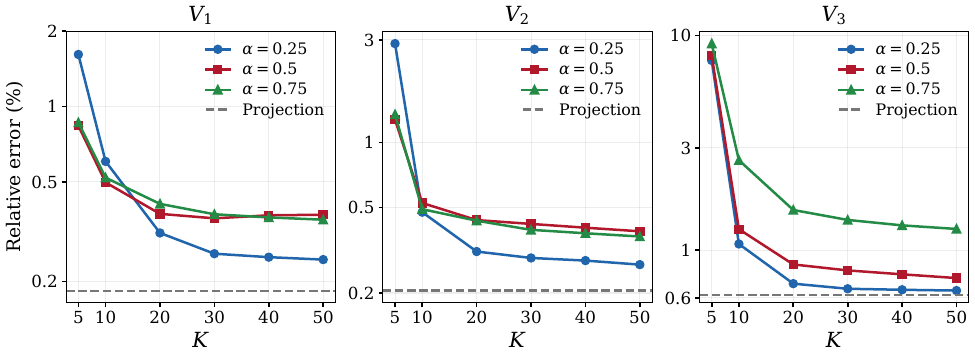}
 \caption{Relative errors versus $K$ on logarithmic vertical axes.
 Dashed lines are the computed least-squares projection errors
 in the $81$-dimensional space on the $161^2$ evaluation grid.}
 \label{fig:dense-bandwidth}
\end{figure}

These examples illustrate the improvement in reconstruction accuracy as the measured frequency band expands.

\FloatBarrier
\endgroup

\section{Conclusion}\label{sec:conclusion}

For every $0<\alpha<1$, we have established high-frequency well-posedness of the fractional Helmholtz equation for real compactly supported $C^1$ potentials.
Under a uniform $C^4$ bound, we have derived increasing stability from complex $L^2$ multifrequency measurements on an enclosing sphere or in an exterior annulus, for all incident directions.
The proof relates the scattered fields to ray integrals through geometric optics, controls the nonlinear phase by frequency combinations and averaging, and extends the frequency range by analytic continuation.

Corresponding estimates under weaker regularity assumptions, with fewer incident directions, or from phaseless measurements remain open.
We hope to report further progress on these problems in future work.

\appendix
\section{Resolvent and analytic continuation estimates}\label{app:resolvent}

This appendix proves results on the direct scattering in Section~\ref{sec:direct} and
Proposition~\ref{prop:extended-domain} in Section~\ref{sec:extended-domain}.

\subsection{Proof of Lemma~\ref{lem:free}}\label{app:free}

\begin{proof}
The fractional-power resolvent identity \cite{MS}, in the form given in \cite[equation (2.22)]{QXZ}, gives
\begin{equation}\label{eq:resolvent-split}
\cR_{\alpha,0}(z)=\frac{z^{2-2\alpha}}{\alpha}R_0(z)+J_\alpha(z).
\end{equation}
The Yukawa kernel of $(\gamma-\Delta)^{-1}$ is $e^{-|x|\sqrt\gamma}/(4\pi|x|)$, giving \eqref{eq:green}--\eqref{eq:J-kernel}.
To estimate the fractional correction and prove its holomorphy, set
\[
D_\alpha(\gamma,z)= (z^{2\alpha}-\gamma^\alpha e^{i\pi\alpha}) (z^{2\alpha}-\gamma^\alpha e^{-i\pi\alpha}).
\]
On a closed subsector of the right half-plane,
\begin{equation}\label{eq:denominator}
|D_\alpha(\gamma,z)|\ge c(\gamma^{2\alpha}+|z|^{4\alpha}).
\end{equation}
Indeed, the argument of $z^{2\alpha}$ stays strictly between $-\pi\alpha$ and $\pi\alpha$.
Scaling $\gamma$ by $|z|^2$ and compactness of the angular interval give the uniform bound.

The Fourier multiplier of $J_\alpha(z)$ is
\[
b_\alpha(q,z)=\frac{\sin(\pi\alpha)}{\pi} \int_0^\infty \frac{\gamma^\alpha}{(\gamma+q^2)D_\alpha(\gamma,z)}\dd\gamma, \qquad q=|\xi|.
\]
Setting $a=|z|$ and $Q=q/a$, \eqref{eq:denominator} yields
\[
|b_\alpha(q,z)|\le Ca^{-2\alpha} \int_0^\infty\frac{s^\alpha}{(s+Q^2)(1+s^{2\alpha})}\dd s \le C(a^2+q^2)^{-\alpha}.
\]
For $Q\le1$, use $s^{\alpha-1}/(1+s^{2\alpha})$ as a majorant.
For $Q>1$, split the integral at $Q^2$; both resulting terms are $O(Q^{-2\alpha})$.
This proves \eqref{eq:J-bound} by Plancherel.
The same bound on compact subsets, together with the scalar Cauchy formula, proves holomorphy in the operator norm $L^2\to H^{2\alpha}$.

The classical localized resolvent is entire and satisfies
\[
\norm{\chi_1R_0(z)\chi_2}_{L^2\to H^l} \le C|z|^{l-1}e^{T(\Im z)_-},\qquad 0\le l\le2,
\]
for $|z|\ge1$; see \cite{DZ}.
Combining this estimate with \eqref{eq:resolvent-split} gives \eqref{eq:free-bound}.

The kernel is locally integrable, with leading singularity $c_\alpha|x|^{2\alpha-3}$ at the origin.
At infinity, its classical part satisfies \eqref{eq:radiation}, while derivatives of the fractional correction have decay $O(|x|^{-3-2\alpha-j})$ at derivative order $j$.
Thus convolution with compactly supported $L^2$ data is outgoing.
For uniqueness, let $w$ be a homogeneous solution in the radiating tempered class with Fourier transform locally integrable near zero.
Then $(|\xi|^{2\alpha}-k^{2\alpha})\widehat w=0$.
Its Fourier support is on $|\xi|=k$.
Multiplication there by the smooth quotient $(|\xi|^2-k^2)/(|\xi|^{2\alpha}-k^{2\alpha})$ gives $(-\Delta-k^2)w=0$.
Classical outgoing uniqueness gives $w=0$.
This agrees with the limiting absorption construction in \cite{ZCV}.
\end{proof}

\subsection{Proof of Theorem~\ref{thm:direct}}\label{app:direct}

\begin{proof}
Choose $\chi\in C_c^\infty$ equal to one near $\supp V$, and consider
\[
T_V(z)=I+\chi\cR_{\alpha,0}(z)V\quad\text{on }L^2(\R^3).
\]
By Lemma~\ref{lem:free}, $T_V(z)-I$ is holomorphic and compact: its range is supported in $\supp\chi$ and gains $2\alpha>0$ derivatives.
If $z$ is in the first quadrant and $T_V(z)v=0$, put $w=-\cR_{\alpha,0}(z)Vv$.
Here the spectral resolvent maps $L^2$ into $H^{2\alpha}$ globally.
Since $\chi w=v$ and $Vw=Vv$, we have
\[
((-\Delta)^\alpha+V-z^{2\alpha})w=0.
\]
Taking the imaginary part of its $L^2$ pairing with $w$ gives $-\Im(z^{2\alpha})\norm w_{L^2}^2=0$, hence $v=w=0$.
The index-zero Fredholm operator is therefore invertible.
Analytic Fredholm theory gives a meromorphic inverse on $\Re z>0$; let $\Sigma_V$ be its positive real singularities.
For a real frequency outside the exceptional set, set
\[
v=T_V(k)^{-1}(\chi e^{ikx\cdot d}),\qquad u=e^{ikx\cdot d}-\cR_{\alpha,0}(k)Vv.
\]
Then $\chi u=v$, proving existence and uniqueness.
The free mapping property gives local $H^{2\alpha}$ regularity and continuous dependence, and defines the localized perturbed resolvent.

Each exceptional frequency corresponds to an $L^2$ eigenfunction.
Suppose $v\ne0$ and $T_V(k)v=0$, and set $w=-\cR_{\alpha,0}(k)Vv$, $f=Vw=Vv$.
The function $f$ is compactly supported and in $L^2$, and
\[
w=-\cR_{\alpha,0}(k)f,\qquad \int_{\R^3}w\overline f\dd x=\int_{\R^3}V|w|^2\dd x\in\R.
\]
The spectral boundary-value formula, paired with the compactly supported $L^2$ function $f$, gives
\begin{equation}\label{eq:shell}
0=\Im\int_{\R^3}(\cR_{\alpha,0}(k)f)\overline f\dd x =\frac{\pi}{(2\pi)^3}\int_{\R^3} \delta(|\xi|^{2\alpha}-k^{2\alpha})|\widehat f(\xi)|^2\dd\xi.
\end{equation}
Compact support makes $\widehat f$ continuously differentiable, so it vanishes on the spectral sphere.
Thus the delta contribution in the outgoing multiplier vanishes, and the quotient $\widehat f(\xi)/(|\xi|^{2\alpha}-k^{2\alpha})$ is bounded in a neighborhood of that sphere.
Away from the sphere the denominator is bounded away from zero on compact sets and has order $2\alpha$ at infinity.
It follows that $w\in H^{2\alpha}(\R^3)$ and
\begin{equation}\label{eq:eigen}
((-\Delta)^\alpha+V)w=k^{2\alpha}w.
\end{equation}

Put $H_0=(-\Delta)^\alpha$, $E=k^{2\alpha}$, and $U_a w(x)=a^{3/2}w(ax)$ for $a>0$.
The identity $H_0U_a=a^{2\alpha}U_aH_0$ and self-adjointness give
\begin{equation}\label{eq:scaling-identity}
(1-a^{2\alpha})E\ip{w}{U_aw} =\ip{w}{(V-a^{2\alpha}V(a\,\cdot))U_aw}.
\end{equation}
Dividing by $1-a^{2\alpha}$ and letting $a\to1$, strong continuity of $U_a$ and uniform convergence of the potential difference quotient yield
\begin{equation}\label{eq:virial}
E\norm w_{L^2}^2= \int_{\R^3}\left(V+\frac{x\cdot\nabla V}{2\alpha}\right)|w|^2\dd x.
\end{equation}
This contradicts $E>E_*(V)$ and proves \eqref{eq:threshold}.
\end{proof}

\subsection{Proof of Lemma~\ref{lem:real-resolvent}}\label{app:real}

\begin{proof}
Set $h=k^{-1}$, $\delta=h^\mu$, and
\[
H_h=|hD|^\mu+\delta V,\qquad A_h=\frac12(x\cdot hD+hD\cdot x),\qquad D=-i\nabla.
\]
Both operators are self-adjoint on their natural domains.
The dilation group preserves $H^\mu(\R^3)$ and gives the exact commutator identity
\begin{equation}\label{eq:dilation-commutator}
[H_h,iA_h] =h\bigl(\mu H_h-\delta(\mu V+x\cdot\nabla V)\bigr).
\end{equation}
For $j=1,2$, put $B_0=H_h$ and $B_j=[B_{j-1},iA_h]$.
Then
\begin{equation}\label{eq:iterated-commutator}
B_j=h^j\bigl(\mu^j|hD|^\mu+\delta(-x\cdot\nabla)^jV\bigr), \qquad \norm{B_j(H_h+i)^{-1}}\le C_jh^j.
\end{equation}
These identities extend from a common core to operators from the domain of $H_h$, with its graph norm, to $L^2$.
They follow from the exact dilation action on the Fourier multiplier, including at zero momentum.

Let $I=[3/4,5/4]$ and denote the spectral projection of $H_h$ onto $I$ by $E_I(H_h)$.
Since $\delta\to0$, \eqref{eq:dilation-commutator} implies, for all sufficiently small $h$,
\[
E_I(H_h)[H_h,iA_h]E_I(H_h)\ge chE_I(H_h),\qquad c>0.
\]
The common core, commutator bounds, and strict lower bound allow us to apply the argument in \cite[Section 4, Lemmas 6--7]{Nakamura}.
We apply the case $\beta=1$, $n=1$, and weight exponent $s=1$, for which (H4) is unnecessary by the remark following Theorem 1 in Section 2 of \cite{Nakamura}.
Its Lemmas 6--7 use only the first two commutator bounds for these parameters.
The passage to spatial weights requires only a Lipschitz bound on $V$.
Consequently, for $\lambda$ in a fixed smaller interval about one,
\begin{equation}\label{eq:dilation-weight}
\norm{\langle A_h\rangle^{-1}(H_h-\lambda-i0)^{-1} \langle A_h\rangle^{-1}}\le Ch^{-1}.
\end{equation}

We obtain spatial weights from \eqref{eq:dilation-weight} by spectral localization.
Choose a real $\varphi\in C_c^\infty(I)$ which is one on a neighborhood of the smaller energy interval.
In the semiclassical Sobolev norm $\norm v_{H_h^s}=\norm{\langle hD\rangle^sv}_2$, we first show that
\begin{equation}\label{eq:spectral-smoothing}
\norm{\varphi(H_h)}_{L^2\to H_h^1}\le C.
\end{equation}
Indeed, multiplication by a uniformly Lipschitz $V$ is bounded on $H_h^s$ for $0\le s\le1$, uniformly for $0<h\le1$, and
\[
\norm v_{H_h^{s+\mu}} \le C_s\big(\norm{|hD|^\mu v}_{H_h^s}+\norm v_{H_h^s}\big), \qquad |hD|^\mu=H_h-\delta V.
\]
Starting from the spectral theorem, iterate this estimate with compactly supported smooth spectral multipliers until the output has one derivative.
The input orders remain below one, where the Lipschitz multiplication bound on $V$ applies.
Functional calculus and \eqref{eq:iterated-commutator} also give $\norm{[A_h,\varphi(H_h)]}\le Ch$.

Put $W=\langle x\rangle^{-1}$, and write $M_f$ for multiplication by $f$.
The ordering identity
\[
A_hW=\sum_{j=1}^3hD_jM_{x_j/\langle x\rangle} +\frac{3ih}{2}M_{1/\langle x\rangle}
\]
and \eqref{eq:spectral-smoothing} imply that $\varphi(H_h)A_hW$ is uniformly bounded.
Commuting $A_h$ past $\varphi(H_h)$ and taking adjoints yield
\[
\norm{\langle A_h\rangle\varphi(H_h)W} +\norm{W\varphi(H_h)\langle A_h\rangle}\le C.
\]
These bounds extend from a core by the closedness of $A_h$.
They allow us to factor
\begin{align*}
W\varphi(H_h)(H_h-\lambda-i0)^{-1}\varphi(H_h)W ={}&\big(W\varphi(H_h)\langle A_h\rangle\big)\\
&\cdot\big(\langle A_h\rangle^{-1}(H_h-\lambda-i0)^{-1} \langle A_h\rangle^{-1}\big) \big(\langle A_h\rangle\varphi(H_h)W\big).
\end{align*}
Its norm is at most $Ch^{-1}$ by \eqref{eq:dilation-weight}.
The complementary term $W(1-\varphi(H_h)^2)(H_h-\lambda-i0)^{-1}W$ is uniformly bounded by the spectral theorem.
We have proved
\[
\norm{W(H_h-\lambda-i0)^{-1}W}\le Ch^{-1}.
\]
At $\lambda=1$, multiplication by $h^\mu$ gives \eqref{eq:weighted-real}.
The constants depend only on a $C^2$ bound.
Mollification and continuity of the localized Fredholm inverse therefore extend the estimate from smooth potentials to the stated class.
\end{proof}

\subsection{Proof of Proposition~\ref{prop:complex-strip}}\label{app:strip}

\begin{proof}
Write the free splitting \eqref{eq:resolvent-split} as
\[
\cR_{\alpha,0}(z)=a(z)R_0(z)+J_z, \qquad a(z)=z^{2-\mu}/\alpha,\quad J_z=J_\alpha(z).
\]
In a fixed closed proper sector of the right half-plane,
\begin{equation}\label{eq:J-derivative}
\norm{J_z}\le C|z|^{-\mu},\qquad \norm{\partial_zJ_z}\le C|z|^{-\mu-1}.
\end{equation}
The derivative estimate follows from \eqref{eq:J-bound} and the Cauchy formula on discs of radius $c|z|$ in a slightly larger sector.
For $\Re z$ large, define
\begin{equation}\label{eq:effective-potential}
B_z=I+J_zV,\quad C_z=I+VJ_z,\quad W_z=C_z^{-1}V=VB_z^{-1},\quad Q_z=a(z)W_z.
\end{equation}
The inverses are uniformly bounded on $L^2$ by a Neumann series.
If $K=\overline{B_r}$, then
\[
W_z=1_KW_z1_K,\qquad \norm{W_z}\le C, \qquad \norm{Q_z}\le C|z|^{2-\mu}.
\]
Here $1_K$ denotes multiplication by the characteristic function of $K$.
The support property follows by expanding $W_z=V-VJ_zV+VJ_zVJ_zV-\cdots$.
At a real frequency $k$, $J_k$ is self-adjoint, so $W_k$ and $Q_k$ are self-adjoint as well.

Set $S_z=(-\Delta-z^2+Q_z)^{-1}_{\rm out}$.
The outgoing resolvents satisfy
\begin{equation}\label{eq:exact-reduction}
\begin{split}
\cR_{\alpha,V}(z)&=B_z^{-1}J_z +a(z)B_z^{-1}S_zC_z^{-1},\\
a(z)S_z&=B_z\cR_{\alpha,V}(z)C_z-J_zC_z.
\end{split}
\end{equation}
To verify them, let $u=\cR_{\alpha,0}(z)(f-Vu)$ and set $w=B_zu-J_zf$.
Direct substitution gives
\[
w=a(z)R_0(z)(C_z^{-1}f-W_zw) =a(z)S_zC_z^{-1}f,
\]
which proves \eqref{eq:exact-reduction} wherever the inverses exist; the identities then extend meromorphically.

For any fixed cutoff $\eta$, the second identity and Lemma~\ref{lem:real-resolvent} imply
\begin{equation}\label{eq:low-classical-real}
\norm{\eta S_k\eta}\le Ck^{-1}.
\end{equation}
Expanding $\eta B_k\cR_{\alpha,V}(k)C_k\eta$ leaves only compactly localized occurrences of $\cR_{\alpha,V}(k)$, since the adjacent factors are either $\eta$ or $V$.
The remaining factors $J_k$ are bounded on global $L^2$.
Division by $a(k)$ changes the power $k^{1-\mu}$ to $k^{-1}$, while the last term in \eqref{eq:exact-reduction} contributes $O(k^{-2})$.

To continue the estimate off the real axis, first freeze the perturbation at $Q_k$.
Write
\[
P^{(k)}=-\Delta+Q_k,\qquad S^{(k)}(z)=(P^{(k)}-z^2)^{-1}_{\rm out}.
\]
This is a self-adjoint bounded perturbation of $-\Delta$ with domain $H^2$.
Its cutoff resolvent is meromorphic in $z$ by analytic Fredholm theory, since the perturbation has input and output supported in $K$ and the localized classical free resolvent gains two derivatives.
At $z=k$, it agrees with $S_k$.

Choose real cutoffs $\rho,\eta,\chi\in C_c^\infty(\R^3)$ with $\rho=1$ near $K$, $\eta=1$ near $\supp\rho$, and $\chi=1$ near $\supp\eta$.
Put
\[
C_\rho=[-\Delta,\rho],\qquad F_z=\eta S^{(k)}(z)\eta,\qquad D_{z,k}=\chi(R_0(z)-R_0(k))\chi.
\]
The two-sided support of $Q_k$ gives $Q_k(1-\rho)=(1-\rho)Q_k=0$.
The cutoff resolvent identity of \cite[Lemma 6.27]{DZ} therefore takes the form
\begin{equation}\label{eq:cutoff-identity}
\begin{split}
F_z-F_k={}&(z^2-k^2)F_z\rho(2-\rho)F_k+\bigl(\eta(1-\rho)+F_zC_\rho\bigr)D_{z,k} \bigl((1-\rho)\eta-C_\rho F_k\bigr).
\end{split}
\end{equation}
The identities
\[
S^{(k)}(z)(1-\rho) =(1-\rho+S^{(k)}(z)C_\rho)R_0(z),
\]
\[
(1-\rho)S^{(k)}(k) =R_0(k)(1-\rho-C_\rho S^{(k)}(k))
\]
follow by applying the respective differential operators.
Substituting $1=\rho(2-\rho)+(1-\rho)^2$ into the ordinary resolvent identity and using these formulas yields \eqref{eq:cutoff-identity}.

In the norms $H_{1/k}^s$, the classical free estimate and Cauchy's formula on a fixed disc give
\begin{equation}\label{eq:free-difference}
\norm{D_{z,k}}_{H_{1/k}^{-r_1}\to H_{1/k}^{r_2}} \le C|z-k|k^{-1},\quad r_1,r_2\in\{0,1\},\quad |z-k|<1;
\end{equation}
see also \cite[equation (6.5.29)]{DZ} after rescaling.
Moreover, $C_\rho$ has norm $O(k)$ both from $H_{1/k}^1$ to $L^2$ and from $L^2$ to $H_{1/k}^{-1}$.
Let $X=\norm{F_z}$ and $Y=\norm{F_k}\le C/k$.
Expanding \eqref{eq:cutoff-identity} and applying the commutators to the free difference give
\[
X\le \frac Ck+C|z-k|\left(X+\frac1k\right).
\]
For the two-commutator term, $\norm{C_\rho D_{z,k}C_\rho}\le Ck|z-k|$, so its contribution is at most $Ck|z-k|XY\le C|z-k|X$.
Thus there is a fixed $b_0>0$ such that
\begin{equation}\label{eq:frozen-strip}
\norm{\eta S^{(k)}(z)\eta}\le Ck^{-1}, \qquad |z-k|<b_0.
\end{equation}
The estimate holds away from the discrete poles and makes each pole in the disc removable.
The constants depend only on the cutoffs and the uniform real-frequency bound in \eqref{eq:low-classical-real}.

We now restore the frequency dependence of the perturbation.
Differentiating \eqref{eq:effective-potential} and using \eqref{eq:J-derivative}, we obtain
\[
\norm{\partial_zW_z}\le Ck^{-\mu-1},\qquad \norm{Q_z-Q_k}\le C|z-k|k^{1-\mu} \quad (|z-k|<b_0).
\]
This difference is supported in $K$ on both sides.
Its product with the localized frozen resolvent has norm at most $C|z-k|k^{-\mu}$ by \eqref{eq:frozen-strip}.
For $k$ large, a Neumann series therefore gives a holomorphic $S_z$ on this disc with the same cutoff bound $C/k$.
Finally, use $B_z^{-1}=I-J_zW_z$ and $C_z^{-1}=I-W_zJ_z$ in \eqref{eq:exact-reduction}.
All occurrences of $S_z$ are compactly localized, and the remaining factors are uniformly bounded on $L^2$.
This gives $\norm{\eta\cR_{\alpha,V}(z)\eta}\le Ck^{1-\mu}$.
The continuations agree on overlapping discs, since they agree in the upper half-plane.
Their union contains \eqref{eq:strip}, after increasing $k_0$ and taking $b<b_0$.
Arbitrary fixed cutoffs follow by choosing $\eta$ to equal one on their supports.
\end{proof}

\subsection{Proof of Corollary~\ref{cor:data-analyticity}}\label{app:data-analyticity}

\begin{proof}
The source $Ve^{izx\cdot d}$ is holomorphic with values in $L^2(B_r\times\Sph)$ and is uniformly bounded in a fixed strip.
By Proposition~\ref{prop:complex-strip} and
\[
u_V^s(\cdot,d,z)=-\cR_{\alpha,V}(z)(Ve^{izx\cdot d}),
\]
the scattered field is holomorphic on every fixed bounded volume observation set, with norm at most $C|z|^{1-\mu}$.
In particular, \eqref{eq:polynomial-data} holds for $\mathcal O=A$, and
\[
\norm{u_V(\cdot,\cdot,z)}_{L^2(B_r\times\Sph)} \le C(1+|z|^{1-\mu}).
\]

For $\mathcal O=\Gamma$, we use the separation $R-r$ to estimate the free kernel directly.
Equations \eqref{eq:J-kernel} and \eqref{eq:denominator} give, for $a>0$,
\[
|J_{\alpha,z}(a)| \le\frac{C|z|^{-2\mu}}a \int_0^\infty\gamma^{\mu/2}e^{-a\sqrt\gamma}\dd\gamma \le C|z|^{-2\mu}a^{-3-\mu}.
\]
Consequently, the operator
\[
T_\Gamma(z)f(x)=\int_{B_r}G_{\alpha,z}(x-y)f(y)\dd y, \qquad x\in\Gamma,
\]
is holomorphic from $L^2(B_r)$ to $L^2(\Gamma)$ and satisfies
\[
\norm{T_\Gamma(z)}\le C(|z|^{2-\mu}+|z|^{-2\mu}).
\]
The kernel is holomorphic and uniformly bounded on the separated compact product for $z$ in compact subsets of the strip, which gives holomorphy in the Hilbert--Schmidt norm.
The identity $u_V^s|_\Gamma=-T_\Gamma(z)(Vu_V)$ proves the assertion, since $2-\mu+(1-\mu)_+<3$.
Finally, the incident and modulating exponentials are bounded holomorphic multipliers on either fixed observation set.
\end{proof}

\subsection{Proof of Proposition~\ref{prop:extended-domain}}\label{app:extended-domain}

Write $H_V=(-\Delta)^\alpha+V$ and $\mu=2\alpha$.
We first establish a spectral estimate for the spherical trace.

\begin{lemma}\label{lem:spectral-H1}
Let $V\in W^{1,\infty}(\R^3;\R)$, let $\lambda\notin\sigma(H_V)$, and put $t=\norm{(H_V-\lambda)^{-1}}_{L^2\to L^2}$.
Then
\begin{equation}\label{eq:spectral-H1}
\norm{(H_V-\lambda)^{-1}f}_{H^1} \le C\left(t\norm f_{H^1} +t^2\norm{\nabla V}_\infty\norm f_{L^2}\right).
\end{equation}
The resolvent is holomorphic on $H^1$ as a function of $\lambda$ in the spectral resolvent set.
\end{lemma}

\begin{proof}
Set $w=(H_V-\lambda)^{-1}f$.
For translation $T_a$ in the $j$th coordinate and $D_a=(T_a-I)/a$, translation invariance of the fractional Laplacian gives
\[
(H_V-\lambda)D_aw=D_af-(D_aV)T_aw.
\]
This identity holds in the natural domain $H^\mu$, which is preserved by translations.
Thus
\[
\norm{D_aw}_{L^2} \le t\norm{\partial_jf}_{L^2} +t^2\norm{\partial_jV}_\infty\norm f_{L^2}.
\]
The difference-quotient characterization of $H^1$ proves \eqref{eq:spectral-H1}.
The resolvent identity and these locally uniform bounds prove holomorphy on $H^1$.
\end{proof}

\begin{proof}[Proof of Proposition~\ref{prop:extended-domain}]
Let $b_0$ be the strip width in Proposition~\ref{prop:complex-strip}.
Choose $0<b<b_0/4$ and $0<\epsilon<\pi/4$ so that $\mu(\pi+\epsilon)<2\pi$.
Increase $L$ so that the domain avoids the chosen cut, $|z|\ge L-b>1$, and
\[
\mu\arg z\in[\mu\pi/2,\mu(\pi+\epsilon)]\Subset(0,2\pi) \quad\hbox{when }z\in\mathcal D_{b,L},\ \Re z\le0.
\]
Since $H_V$ is self-adjoint and $\sigma(H_V)\subset[-F,\infty)$, a further increase of $L$ gives
\begin{equation}\label{eq:left-spectral-bound}
\norm{(H_V-z^\mu)^{-1}}_{L^2\to L^2} \le C|z|^{-\mu},\qquad \Re z\le0, \quad z\in\mathcal D_{b,L}.
\end{equation}
Indeed, the distance of $e^{i\theta}$ from $[0,\infty)$ has a positive minimum on the displayed angular interval, and
\[
\operatorname{dist}(z^\mu,[-F,\infty)) \ge c_\mu|z|^\mu-F\ge\tfrac12c_\mu|z|^\mu.
\]

In the first quadrant, $\sin(\mu\theta)\ge c_\mu\sin\theta$ for $0<\theta\le\pi/2$, and hence
\begin{equation}\label{eq:upper-spectral-bound}
\norm{(H_V-z^\mu)^{-1}}_{L^2\to L^2} \le\frac{1}{\Im(z^\mu)} \le\frac{C|z|^{1-\mu}}{\Im z}.
\end{equation}
Use the spectral resolvent on the two open sets $\mathcal D_{b,L}\cap\{\Im z>0\}$ and $\mathcal D_{b,L}\cap\{\Re z<0\}$, and the outgoing resolvent on
\[
\mathcal D_{b,L}\cap\{\Re z>k_c,\ \Im z<2b\},
\]
where $k_c$ is the lower endpoint of the fixed strip.
Taking $L>\sqrt{k_c^2+9b^2}+1$ makes these sets an open cover of $\mathcal D_{b,L}$.
The spectral and outgoing resolvents agree on their first-quadrant overlap and therefore define a holomorphic extension.
The strip estimate, \eqref{eq:left-spectral-bound}, and \eqref{eq:upper-spectral-bound} for $\Im z\ge b$ prove \eqref{eq:extended-resolvent-bound}.

On the spectral pieces, write
\[
u_V^s(z)=-\bigl(H_V-z^\mu\bigr)^{-1}f_V(z), \qquad f_V(z,x,d)=V(x)e^{izx\cdot d}.
\]
The source is entire with values in $H_x^1L_d^2$ and satisfies
\[
\norm{f_V(z)}_{L^2_{x,d}}\le Ce^{r|\Im z|},\qquad \norm{f_V(z)}_{H_x^1L_d^2} \le C(1+|z|)e^{r|\Im z|}.
\]
On $\Re z\le0$ and on $\Re z\ge0$, $\Im z\ge b$, the spectral resolvent norm is at most $C(1+|z|)$.
Lemma~\ref{lem:spectral-H1} consequently yields
\begin{equation}\label{eq:global-spectral-field-H1}
\norm{u_V^s(z)}_{H_x^1L_d^2} \le C(1+|z|)^2 e^{r(\Im z)_+}.
\end{equation}
Here $e^{r|\Im z|}\le e^{rb}e^{r(\Im z)_+}$ in the domain.
The bounded trace map $H^1(\R^3)\to L^2(\partial B_R)$ gives the spherical bound, and the $L^2$ restriction map gives the annular bound.
These bounded maps and the $H^1$ resolvent holomorphy also give $Y$-valued holomorphy.

On the first-quadrant overlap, the compact-source free-kernel formula for this field is smooth at the receivers, since $R-r>0$.
Its $H^1$ trace therefore equals the spherical restriction in Corollary~\ref{cor:data-analyticity}.
Hence the spectral traces and the outgoing receiver fields define one holomorphic continuation of the measured field.
On the remaining bounded-height strip, Corollary~\ref{cor:data-analyticity} supplies the cubic bound.
Combining it with \eqref{eq:global-spectral-field-H1} proves \eqref{eq:extended-data-bound}.
\end{proof}

\section{Symbol and free propagation estimates}\label{app:geometric-optics}

We prove the finite symbol expansions and free propagation bounds used in Section~\ref{sec:inverse}.
Throughout this appendix, $\mu\in(0,2)$ is fixed.

\subsection{Finite expansions of the smooth symbol}\label{app:finite-kernel}

Use the symbol $p_{\rm sm}$ chosen in Section~\ref{sec:phase-construction}, and let $K=\mathcal F^{-1}p_{\rm sm}$.
Away from zero, its kernel satisfies
\begin{equation}\label{eq:smooth-kernel-bounds}
|K(v)|\le C|v|^{-3-\mu}\quad(0<|v|\le1),\qquad |K(v)|\le C_M|v|^{-M}\quad(|v|\ge1)
\end{equation}
for every sufficiently large integer $M$.
The first bound follows from the finite-part kernel of $|D|^\mu$, since $|\xi|^\mu-p_{\rm sm}$ is compactly supported and integrable.
The second follows by integration by parts in frequency, since sufficiently high derivatives of $p_{\rm sm}$ are integrable.
The finite-part subtraction at zero involves only the value and first derivatives of the function on which the kernel acts.
For a fixed radial cutoff $\vartheta$ equal to one near zero,
\[
\langle K,F\rangle=c_*F(0)+\int K(v)\bigl(F(v)-\vartheta(v)F(0)-\vartheta(v)v\cdot\nabla F(0)\bigr)\dd v.
\]
Here $c_*$ depends only on the fixed symbol and cutoff; radial symmetry removes first-order terms supported at zero, and the order $\mu<2$ excludes higher-order ones.
On a function vanishing to second order at zero, it acts by ordinary integration against $K$.

Let $S$ be a real phase with bounded momentum $\Xi=\nabla S$, and let $a$ have bounded derivatives of the indicated orders.
Set
\[
F_t(x,v)=e^{i\Phi_t(x,v)}a(x-tv),\qquad \Phi_t(x,v)=-\int_0^1v\cdot\Xi(x-stv)\dd s.
\]
Then
\[
e^{-iS(x)/h}p_{\rm sm}(hD)(e^{iS/h}a)(x)=\langle K,F_h(x,\cdot)\rangle.
\]
If $a,\Xi\in C^L$ and $L\ge2$, direct differentiation gives
\[
|\partial_t^LF_t(x,v)|\le C|v|^L(1+|v|^L).
\]
Taylor's formula through degree $L-1$ has a remainder $\mathcal E_{L,h}$ bounded by $Ch^L|v|^L(1+|v|^L)$.
Moreover,
\[
F_t(x,0)=a(x),\qquad \nabla_vF_t(x,0)=-i\Xi(x)a(x)-t\nabla a(x),
\]
so $\mathcal E_{L,h}$ and its first $v$ derivatives vanish at zero.
Since $L\ge2>\mu$, \eqref{eq:smooth-kernel-bounds} yields
\begin{equation}\label{eq:kernel-taylor-budget}
|\langle K,\mathcal E_{L,h}\rangle|\le Ch^L\int|K(v)|\,|v|^L(1+|v|^L)\dd v\le Ch^L.
\end{equation}
Thus
\begin{equation}\label{eq:finite-symbolic-remainder}
e^{-iS/h}p_{\rm sm}(hD)(e^{iS/h}a)=\sum_{\ell=0}^{L-1}h^\ell\mathscr L_{\ell,S}a+h^Lr_{L,h},\qquad \norm{r_{L,h}}_\infty\le C.
\end{equation}
The expansion uses $L$ derivatives of $a$ and $\Xi$; differentiation in $x$ through order $q$ requires $L+q$ derivatives.
If $a=O_{C^L}(\delta)$, the remainder is $O(\delta)$.
The same conclusion holds when $a-1$ and $\Xi-d$ are $O_{C^L}(\delta)$, since every positive $t$ derivative then contains a differentiated amplitude or a positive spatial derivative of $\Xi$.

The bilinear phase representation also gives the coefficients.
With $A_x(z)=\int_0^1\Xi(x+tz)\dd t$, they are
\begin{equation}\label{eq:coefficient-formula}
\mathscr L_{\ell,S}a(x)=\frac{1}{i^\ell}\sum_{|\beta|=\ell}\frac1{\beta!}\left.\partial_z^\beta\bigl((\partial_\xi^\beta p_{\rm sm})(A_x(z))a(x+z)\bigr)\right|_{z=0}.
\end{equation}
The change of variables $\eta=\xi-A_x(z)$ reduces the phase to $-z\cdot\eta$, whose diagonal Taylor coefficients give \eqref{eq:coefficient-formula}.
The coefficient of a derivative of order $q$ of $a$ uses derivatives of $\Xi$ up to order $\ell-q$.
In particular,
\[
\mathscr L_{0,S}a=p_{\rm sm}(\Xi)a,\qquad \mathscr L_{1,S}a=i^{-1}\mathcal T_Sa,\qquad \mathcal T_Sa=v_S\cdot\nabla a+\tfrac12(\operatorname{div}v_S)a,
\]
where $v_S=\nabla_\xi p_{\rm sm}(\Xi)$.

For a fixed smooth cutoff $\chi$, the commutator has the representation
\[
e^{-iS/h}[p_{\rm sm}(hD),\chi](e^{iS/h}a)=\langle K,e^{i\Phi_h}(\chi(x-hv)-\chi(x))a(x-hv)\rangle.
\]
The identity
\[
\chi(x-hv)-\chi(x)=-h\int_0^1v\cdot\nabla\chi(x-thv)\dd t
\]
supplies one power of $h$ before Taylor expansion.
Consequently, for $L\ge2$,
\begin{equation}\label{eq:commutator-remainder}
e^{-iS/h}[p_{\rm sm}(hD),\chi](e^{iS/h}a)=\sum_{\ell=1}^{L-1}h^\ell[\mathscr L_{\ell,S},\chi]a+h^Lr^\chi_{L,h}.
\end{equation}
The $C^q$ bound for $r^\chi_{L,h}$ requires $L-1+q$ derivatives of $a$ and $\Xi$.
The Taylor remainder vanishes to order $L$ in $v$, making its kernel integral convergent at zero.
If $a=O_{C^{L-1+q}}(\delta)$, this remainder is $O_{C^q}(\delta)$.

These expansions hold uniformly on the fixed construction sets.
The phase may be extended by a fixed cutoff applied to $S-x\cdot d$, and amplitudes by cutoffs applied to their deviations from their incoming values.
For small $\delta$, the extensions retain $|\nabla S|\ge1/2$ and the stated finite derivative bounds.
Outside a fixed neighborhood of the spatial support of an input, \eqref{eq:smooth-kernel-bounds} gives arbitrarily small powers of $h$ with any fixed polynomial spatial weight.
The same holds for a commutator outside a fixed neighborhood of the cutoff derivative support, where its kernel is evaluated at a fixed positive separation.

\subsection{Separated source and receiver estimates}\label{app:receiver-estimates}

We use the decomposition
\[
R_{h,0}=\frac2\mu(-h^2\Delta-1-i0)^{-1}+j_\mu(hD),\qquad j_\mu(\xi)=\frac1{|\xi|^\mu-1}-\frac{2/\mu}{|\xi|^2-1}.
\]
The convolution kernel $J_h$ of the second term is positive, and the Yukawa representation gives
\begin{equation}\label{eq:surface-j-mass}
\norm{J_h}_1=\frac2\mu-1,\qquad |\partial^\beta J_h(z)|\le C_\beta h^\mu|z|^{-3-\mu-|\beta|}\quad(|z|\ge c>0).
\end{equation}

For sources in the potential support,
\begin{equation}\label{eq:separated-sphere-L2}
\norm{R_{h,0}}_{L^2(B_r)\to L^2(\Gamma)}\le Ch^{-1}.
\end{equation}
Parameterize a receiver patch by $X(x')$, $x'\in\R^2$, and put $\varphi(x',y)=|X(x')-y|$.
Fix $r<r_1<R$, extend the source by zero, and insert a smooth source cutoff equal to one on $B_r$ and supported in $B_{r_1}$.
With $\omega=(X-y)/|X-y|$, its mixed derivative matrix is
\[
\partial_{x'}\partial_y\varphi=-|X-y|^{-1}(D_{x'}X)^T(I-\omega\otimes\omega).
\]
This matrix has rank two, since $\omega\cdot n_X\ge(R-r_1)/(R+r_1)>0$ on the cutoff support.
After a finite partition, write $y=(y',t)$ so that the corresponding two by two minor is uniformly invertible.
For fixed $t$, let $T_{h,t}$ have kernel $e^{i\varphi(x',y',t)/h}a(x',y',t)$, where $a$ is the smooth compactly supported geometric amplitude.
On sufficiently small patches,
\[
|\nabla_{y'}(\varphi(x',y',t)-\varphi(z',y',t))|\ge c|x'-z'|.
\]
Integration by parts bounds the kernel of $T_{h,t}T_{h,t}^*$ by $C_N(1+|x'-z'|/h)^{-N}$ for any fixed $N$.
Schur's test in two dimensions yields $\norm{T_{h,t}}\le Ch$, uniformly in $t$.
Integration in $t$ and Cauchy--Schwarz give the same bound from the source volume to the receiver patch.
The factor $h^{-2}$ in the classical kernel proves its contribution to \eqref{eq:separated-sphere-L2}.
The separated bound in \eqref{eq:surface-j-mass} controls the remaining term by $Ch^\mu$.

The same estimate holds with any fixed source radius smaller than $R$.

\subsection{Proof of Lemma~\ref{lem:directed-free-trace}}\label{app:directed-trace}

\begin{proof}
Rotate coordinates so that $d=e_3$, and write $\xi=(\eta,\tau)$, with $\eta\in\R^2$.
On the support of $q_d$, put
\[
a_+(\eta)=\sqrt{1-|\eta|^2},\qquad
b(\eta,\tau)=\frac{(|\eta|^2+\tau^2)^{\mu/2}-1}{\tau-a_+(\eta)}.
\]
The quotient extends smoothly across $\tau=a_+(\eta)$ and is positive and bounded away from zero on a fixed neighborhood of $(0,1)$.
Define $m=q_d/b$ on this neighborhood and extend it by zero outside the support of $q_d$.
The outgoing multiplier factors as
\begin{equation}\label{eq:directed-factorization}
R_{h,0}q_d(hD)=(hD_t-a_+(hD_p)-i0)^{-1}m(hD).
\end{equation}
Here $a_+$ may be extended to a smooth real function outside a fixed small neighborhood of the origin.
Choose $\psi\in C_c^\infty(\R^2)$ equal to one on the transverse frequency projection of $\operatorname{supp}m$, and supported in $|\eta|<2\varepsilon$.
For $G=m(hD)f$, the retarded formula is
\begin{equation}\label{eq:directed-retarded}
R_{h,0}q_d(hD)f(p,t)
=\frac{i}{h}\int_{-\infty}^t
e^{i(t-s)a_+(hD_p)/h}\psi(hD_p)G(\cdot,s)(p)\dd s.
\end{equation}
One-dimensional outgoing Fourier inversion in $\tau$ proves the identity for smooth compactly supported $f$; the estimates below extend it to the stated class.

Choose fixed radii
\[
a<a_1<b<c<R.
\]
The portion of $\Gamma$ with $|p|<c$ consists of the two graphs
\[
t=\rho_\pm(p)=\pm\sqrt{R^2-|p|^2}.
\]
Let $\chi(p)$ be a smooth cutoff supported in $|p|<c$ and equal to one on $|p|\le b$.
For either graph, consider the variable-time propagator
\[
T_\rho v(p)=\chi(p)(2\pi h)^{-2}
\int e^{i((p-z)\cdot\eta+\rho(p)a_+(\eta))/h}
\psi(\eta)v(z)\dd z\dd\eta.
\]
Its phase in $(p,\eta)$ has mixed derivative matrix
\[
I+\nabla\rho(p)\otimes\nabla a_+(\eta).
\]
Since $\nabla\rho$ is bounded on $|p|\le c$, a sufficiently small $\varepsilon$ keeps the determinant uniformly away from zero.
The semiclassical Fourier transform is unitary, and the two-dimensional nondegenerate oscillatory integral estimate therefore gives
\begin{equation}\label{eq:variable-time-L2}
\|T_\rho\|_{L^2(\R^2)\to L^2(\R^2)}\le C.
\end{equation}
After finite localization in $(p,\eta)$, the kernel of the unnormalized oscillatory integral operator composed with its adjoint is bounded by $C_M(1+|p-p'|/h)^{-M}$.
Integration by parts in $\eta$ gives this bound, since the gradient of the phase difference is at least $c|p-p'|$.
Schur's test gives an $O(h)$ norm for the unnormalized two-dimensional operator, and its Fourier normalization cancels this factor.

Inserting $t=\rho(p)$ in \eqref{eq:directed-retarded} gives
\[
\chi(p)R_{h,0}q_d(hD)f(p,\rho(p))
=\frac{i}{h}\int_{\R}
\mathbf 1_{s<\rho(p)}T_\rho e^{-is a_+(hD_p)/h}G(\cdot,s)(p)\dd s.
\]
For each $s$, the characteristic function multiplies the output and has $L^2$ operator norm at most one.
The fixed-time propagator is unitary.
Minkowski's inequality and \eqref{eq:variable-time-L2} give
\[
\|\chi R_{h,0}q_d(hD)f|_{t=\rho(p)}\|_2
\le Ch^{-1}\int_\R\|G(\cdot,s)\|_{L^2_p}\dd s.
\]
The kernel $M_h$ of $m(hD)$ has a uniform $L^1(\R^3)$ bound, so mixed-norm Young's inequality yields
\[
\|G\|_{L^1_sL^2_p}\le C\|f\|_{L^1_sL^2_p}\le C_T\|f\|_2.
\]
The graph surface measure is uniformly comparable with $dp$, proving the estimate on both caps.

On the belt $|p|\ge b$, choose $T_1>T$ and a fixed smooth cutoff $\chi_1$ equal to one on a neighborhood of $K_d$, with
\[
\operatorname{supp}\chi_1\subset\{|p|<a_1,\ |s|<T_1\}.
\]
Write $G=G_{\rm near}+G_{\rm far}$, where $G_{\rm near}=\chi_1G$.
Since $m$ is smooth and compactly supported, its kernel is Schwartz, and the positive distance between $K_d$ and the support of $1-\chi_1$ gives, for every $M$,
\begin{equation}\label{eq:directed-filter-tail}
\|G_{\rm far}\|_{L^1(\R^3)}\le C_Mh^M\|f\|_2.
\end{equation}
Split \eqref{eq:directed-retarded} accordingly, retaining $\psi(hD_p)$ on both terms.
The propagator kernel with this cutoff is bounded in absolute value by $Ch^{-2}$ for all time differences.
Thus the far contribution to \eqref{eq:directed-retarded} is bounded pointwise by $C_Mh^{M-3}\|f\|_2$.

For the near contribution, $|z|<a_1$ and $|s|<T_1$, while $x\in\Gamma$ gives $|t|\le R$.
The transverse phase gradient in the propagator kernel is
\[
p-z+(t-s)\nabla a_+(\eta).
\]
Reduce $\varepsilon$ so that, uniformly on the support of $\psi$,
\[
(R+T_1)|\nabla a_+(\eta)|\le\tfrac12(b-a_1).
\]
Thus the gradient has magnitude at least $(b-a_1)/2$ on the belt.
Repeated integration by parts in $\eta$ bounds this kernel by $C_Mh^{M-2}$ for any prescribed order.
The time-ordering factor $\mathbf 1_{s<t}$ is independent of $\eta$.
Since $\|G_{\rm near}\|_1\le C\|G\|_2\le C\|f\|_2$, the near contribution is at most $C_Mh^{M-3}\|f\|_2$.
Together with \eqref{eq:directed-filter-tail}, this gives arbitrarily small powers of $h$ on the belt and proves \eqref{eq:microlocal-sphere} on the whole sphere.

Finally suppose the source is downstream of the observation set as in \eqref{eq:backward-cone}.
For every $x\in K_O$, the retarded formula uses only $s<t\le T_-$, which is separated from the source support $s\ge T_+$ by a fixed positive distance.
Schwartz kernel decay therefore gives
\[
\|\mathbf 1_{s\le T_-}G\|_{L^1(\R^3)}\le C_Mh^M\|f\|_2.
\]
The absolute $Ch^{-2}$ bound for the transversely localized propagator and the factor $h^{-1}$ in \eqref{eq:directed-retarded} prove \eqref{eq:backward-cone}.
The constants are uniform in $d$ by rotation, since the geometric parameters and symbol bounds are fixed.
\end{proof}

\section{A continuation estimate}\label{app:continuation}

The estimate below uses a fixed initial subinterval of the measurements, so its constants are independent of the upper measured frequency.

\begin{lemma}\label{lem:exterior-continuation}
Let $b>0$, $L>b$, and $K_0\ge2L+2$.
Suppose that $F$ is holomorphic with values in a complex Banach space $\mathcal H$ on
\[
\mathcal D_{b,L}=\{z:\Im z>-b,\ |z+ib|>L\}, \qquad \norm{F(z)}_{\mathcal H}\le1.
\]
If $\norm{F(k)}_{\mathcal H}\le\eps$ on $[K_0,K_0+1/2]$, where $0<\eps<1$, then
\begin{equation}\label{eq:exterior-continuation}
\norm{F(q)}_{\mathcal H} \le\exp\left(-\frac{c|\log\eps|}{(1+q)^2}\right), \qquad q\ge K_0+1.
\end{equation}
Here $c>0$ depends only on $b,L,K_0$.
\end{lemma}

\begin{proof}
The Joukowski map
\[
J(z)=\zeta+\frac{L^2}{\zeta},\qquad\zeta=z+ib,
\]
maps $\mathcal D_{b,L}$ conformally onto the upper half-plane $\C_+$.
Indeed, $\Im J(z)=\Im\zeta(1-L^2/|\zeta|^2)>0$.
The roots of $\zeta^2-J(z)\zeta+L^2=0$ have product $L^2$, and the unique root outside $|\zeta|=L$ lies in the upper half-plane.
This defines the holomorphic inverse.

Set $I=[K_0,K_0+1/2]$, $\eta(t)=J(t)$, and $\mathcal A=\eta(I)$.
On this compact smooth arc,
\begin{equation}\label{eq:continuation-arc}
\frac{3b}{4}\le\Im\eta(t)\le b, \qquad \frac{\mathrm d}{\mathrm dt}\Re\eta(t) =1-L^2\frac{t^2-b^2}{(t^2+b^2)^2}\ge\frac34.
\end{equation}
Let
\[
g(w,\eta)=\frac1{2\pi} \log\left|\frac{w-\overline\eta}{w-\eta}\right| =\frac1{4\pi}\log\left(1+ \frac{4\Im w\Im\eta}{|w-\eta|^2}\right), \qquad U(w)=\int_I g(w,\eta(t))\dd t.
\]
The positive potential $U$ is harmonic off $\mathcal A$, bounded and continuous on $\C_+$, zero on the real boundary, and tends to zero at infinity.
The monotonicity in \eqref{eq:continuation-arc} bounds the local logarithmic singularity by $C+|\log|t-t_0||$ with integrals over $|t-t_0|<h$ bounded by $Ch(1+|\log h|)$, uniformly in $w$ near the arc.
Here $t_0$ is the inverse image of the projection of $\Re w$ onto the interval $\Re\eta(I)$.
Splitting the integral at $h$ proves continuity and boundedness.
The Green function gives the boundary limits and decay.

Choose a fixed $B\ge\max(1,\sup U)$ and write $u=U/B$.
For $\ell\in\mathcal H^*$ with $\norm\ell\le1$, put
\[
v(w)=\log|\ell(F(J^{-1}(w)))|, \qquad E=|\log\eps|.
\]
This function is subharmonic, bounded above by zero, and at most $-E$ on $\mathcal A$.
Since $0\le u\le1$, the function $v+Eu$ has upper boundary limit at most zero on the arc, the real axis, and at infinity.
The maximum principle on $\C_+\setminus\mathcal A$ yields $v\le-Eu$.
On bounded exhaustions, the added circular boundary contributes at most $E\sup_{|w|=T}u(w)\to0$.

For $q\ge K_0+1$, the point $w_q=J(q)$ satisfies $3b/4\le\Im w_q\le b$, while $|w_q-\eta(t)|^2\le C(1+q)^2$ for $t\in I$.
The Green function therefore gives
\[
U(w_q)\ge\frac{c_1}{(1+q)^2}.
\]
It follows that $|\ell(F(q))|\le\exp(-cE/(1+q)^2)$.
Taking the supremum over such $\ell$ proves \eqref{eq:exterior-continuation}.
\end{proof}

\section*{Acknowledgements}
X. Xu is partially supported by the National Key Research and Development Program of China (No. 2024YFA1012303), National Natural Science Foundation of China (No. 12525112), and the Open Research Project of Innovation Center of Yangtze River Delta, Zhejiang University. Y. Zhao is supported by National Natural Science Foundation of China (No. 12371421).

\section*{Conflict of Interest Statement}
The authors declare that they have no conflict of interest.

\section*{Data Availability Statement}
The numerical experiments use synthetic data generated using the models and parameters described in this article.
The data supporting the numerical results are available from the corresponding author upon reasonable request.

\end{document}